\documentclass[final]{siamltex}
\usepackage{amsfonts}
\usepackage{amssymb}

\usepackage{titlesec}
\titleformat{\section}
  {\normalfont\Large\bfseries\centering}{\thesection}{1em}{}

\title{Rational curves on hypersurfaces}

\author{ B. Wang\\
Aug 20, 2014}

\begin{document}

\maketitle

\begin{abstract} In this paper, we prove three related results;\par

(1) Extension of our result in [10] to all generic  hypersurfaces. More precisely, 
the normal sheaf of a generic rational map $c_0$  to a generic hypersurface $X_0$ of $\mathbf P^n, n\geq 4$ has
a vanishing higher cohomology, 
\begin{equation}
H^1(N_{c_0/X_0})=0.
\end{equation}
 \par

As applications we give
 \par

(2) A solution to a Voisin's conjecture [9] on a covering of  a generic hypersurface by rational curves  \par
(3) A classification of rational curves on hypersurfaces of general type--a solution  to another  
Voisin's conjecture [9]. 

\end{abstract}

\section{Introduction}\quad
\smallskip

\subsection{Statement}\smallskip\quad

 We work over complex numbers, $\mathbb C$. A general or generic  hypersurface is referred  to a hypersurface as a point in a complement of
a countable union of  proper closed subsets of the space of all hypersurfaces. 
Let
$\mathcal L=\mathcal O_{\mathbf P^n}(1)$, and 
$$X_0=div(f_0)$$
where $f_0\in H^0(\mathcal L^h)$ is generic. 
  Let 
$$ c_0: \mathbf P^1 \to C_0\subset  X_0\subset \mathbf P^n$$ be birational  onto
an irreducible rational curve $C_0$. 
Let $N_{c_0/X_0}$ be the normal sheaf of  the birational morphism $c_0$.\footnote{Throughout the paper, one should be cautious about a fact that
$c_0^\ast(T_{X_0}) $ is a locally free sheaf, but $N_{c_0/X_0}$ may not be. Thus the notation $N_{c_0/X_0}$ is always refereed to a
sheaf module, the normal sheaf of the morphism.}
 Recall that it is a sheaf uniquely determined by the exact sequence
of sheaves over $\mathbf P^1$, 
$$\begin{array}{ccccccccc}
0& \rightarrow & T_{\mathbf P^1}&\stackrel{(c_0)_\ast} \rightarrow & c_0^\ast(T_{X_0}) &\rightarrow & 
N_{c_0/X_0} &\rightarrow &0, \end{array}$$
where $(c_0)_\ast$, the differential of $c_0$ is an injective morphism of sheaf modules.

\bigskip

\begin{theorem} (Main Theorem).  Assume $ n\geq 4$.  Also assume that $X_0$ is a generic hypersurface containing $C_0$, and
$c_0$ is generic in an irreducible component of $Hom(\mathbf P^1, X_0)$.   
 Then
\begin{equation} H^1(N_{c_0/X_0})=0.\end{equation}

 \end{theorem}
\bigskip

\subsection{Applications}\smallskip\quad

The main application of the result  is in the area of rational curves on hypersurfaces. 
General hypersurfaces $X_0$ of degree $h$ in $\mathbf P^n$ can be divided  into three different classes:\par

(1) Fano variety,  in that case $h<n+1$ and there are lots of rational curves of each degree,  \par
(2) Calabi-Yau variety,  in that case $h=n+1$ and
there are ``fewer" rational curves of each degree,  \par
(3) a variety of general type,  in that
case $h> n+1$ and there are no rational curves for large $h$. 
\par

A general question:  what are the structures  of the Hilbert schemes Hilb$(X_0)$ of  rational curves on $X_0$? 

 \par

 For the class (1), there are results  of J. Harris, M. Roth and J. Starr (HRS)  in [4](using different techniques from us),  and
  a conjecture about Hilbert scheme Hilb$_d(X_0)$ of rational curves of degree $d$ by Coskun, Harris and Starr, who extend the results in [4] to almost all Fano hypersurfaces. For the classes (2), (3), the situations are quite different. However 
in a view of deformation theory,  $ H^1(N_{c_0/X_0})$ is the obstruction space of the deformation of $c_0$ on $X_0$ if the rational map $c_0$ is smooth.  So we try to identify the cohomology $ H^1(N_{c_0/X_0})$ as a unified invariant. In this view 
the result, \begin{equation} H^1(N_{c_0/X_0})=0.\end{equation}
has many implications. For examples if $c_0$ is a smooth embedding, 
it implies that \par
(1) deformation of $c_0$ on $X_0$ has no obstruction,\par
(2) deformation of $c_0$  over the complex moduli of $X_0$  has no obstruction.
 \par
These are immediate consequences of theorem 1.1.  They all assert  Hilbert schemes above have expected dimensions. More directly 
propositions 1.2, 1.3 below (or theorem 1,1 ), in the case $n\geq 4$, can be used  to reproduce almost all results in [4], and furthermore to prove the conjecture made by Coskun, Harris and Starr \footnote{It can't be directly used to prove the irreducibility of Hilb$_d(X)$ for $X$ in class (1) (as proved and conjectured).
Besides,  the irreducibility in [4] clearly can not extend to classes (2), (3), while all the other results can.}.  This will address the question for class (1). 
 However the detailed discussion of this will be given elsewhere.  In this paper we'll only concentrate on the classes (2), (3). 
 Most of the known work and conjectures in this part (Classes (2), (3)) were nicely summarized by Voisin in [9].
In [9], Voisin studied the rational curves and elliptic curves on a variety $X_0$.  As consequences she formulated many 
interesting conjectures. Here we list two of them.\par
 (I) In a general Calabi-Yau hypersurface of dimension $\geq 3$, rational 
curves cover a countable union of Zariski closed proper algebraic subsets of codimension $\geq 2$. \par
(II)  If a generic hypersurface  is of general type,  
the degrees of rational curves on it are bounded.\par

In this paper, we apply theorem 1.1 to show that Voisin's conjectures (I), (II)  are correct. 

\subsection{Outline of the proof of Main theorem}\smallskip\quad

There is an easy reduction (see section 2.2) that shows that  it suffices to prove theorem 1.1 for  the Calabi-Yau hypersurface $X_0$, i.e.
$$n+1=h.$$
Let's outline the proof for this Calabi-Yau's case.  
Let $S=\mathbf P(H^0(\mathcal O_{\mathbf P^n}(h)))$ be the space of all hypersurfaces of degree $h$. 
Theorem 1.1 is stated in terms of rational curves and hypersurfaces in projective space. 
But in this paper we'll stick with the affine space for the simplicity. So let
$$\mathbb C^{(n+1)(d+1)}$$ be the
vector space, 
$$( H^0(\mathcal O_{\mathbf P^1}(d))^{\oplus n+1}$$
whose open subset parametrizes  the set of maps $$\mathbf P^1\to \mathbf P^n$$ whose push-forward cycles have degree $d$.\footnote {The automorphism of $\mathbf P^1$ induces a
$PGL(2)$ group action on $\mathbf P(\mathbb C^{(n+1)(d+1)})$.  Let $$PGL(2)( c_0)\subset \mathbf P(\mathbb C^{h(d+1)})$$ be the orbit of $c_0\in \mathbf P(\mathbb C^{h(d+1)})$. } 
Throughout the paper, we let $$M= \mathbb C^{(n+1)(d+1)}.$$ 
Let $M_d$ be the subset that consists of all birational-to-its-image maps $c$ whose push-forward cycles $c_\ast([\mathbf P^1])$ have degree $d$. 
$M$ has affine coordinates.  Assume $c_0^\ast(f_0)=0$ for $c_0\in M_{d}$ as in theorem 1.1.
Let $\mathbb L\subset S$ be an open set of the plane spanned by hypersurfaces $f_0, f_1, f_2$, where $f_0, f_1, f_2$ are generic in $S$.
Let  $$\Gamma_{\mathbb L}\ni (c_0, [f_0])$$ be an irreducible component of the incidence scheme
\begin{equation}\{(c, [f])\subset M\times \mathbb L: c^\ast(f)=0\}
\end{equation} 
that is onto $\mathbb L$.  We assume $\Gamma_{\mathbb L}$ exists. Let $P$ be the projection $$\Gamma_{\mathbb L}\to M.$$
The idea of the proof is to show that the scheme,
$$P(\Gamma_{\mathbb L})\subset M$$ is
a reduced, irreducible quasi-affine scheme of dimension $h+1$.  The method is straightforward to show its defining polynomials at a generic point have non-degenerate Jacobian matrix (by that we mean it has full rank).  See
definition 2.6 below for the precise definition of a Jacobian matrix.   All differentials and partial derivatives used throughout the paper are in algebraic sense (because all functions used are holomorphic).  In the following we describe its defining polynomials and a differential form representing
the Jacobian matrix. 

Choose generic $hd+1$ distinct points  $t_i\in \mathbf P^1$ (generic in $Sym^{hd+1}(\mathbf P^1)$). In the following (1.4), we use
$t_i$ to denote a nonzero point in $\mathbb C^2$ with two coordinates $( \sigma_1^i, \sigma_2^i)$, i.e $P(\mathbb C^2)=\mathbf P^1$.
Next we consider differential 1-forms $\phi_i$ on $M$:
  
\begin{equation} \phi_i= d
\left|  \begin{array}{ccc} f_2(c(t_i)) & f_1(c(t_i)) & f_0(c(t_i))\\
f_2(c(t_1)) & f_1(c(t_1)) & f_0(c(t_1))\\
f_2(c(t_2)) & f_1(c(t_2)) & f_0(c(t_2))
\end{array}\right|\end{equation}

for $i=3, \cdots, hd+1$, and variable $c\in M$, where $|\cdot |$ denotes the determinant of a 
matrix.  Notice $\phi_i$ are uniquely defined provided  the quintics $f_i$ are in an affine open set of $S$ ( $\phi_i$ is not invariant under the $GL(2)$ action of $\mathbb C^2$ ).  Let \begin{equation} \omega=\wedge_{i=3}^{hd+1}\phi_i \in H^0(\Omega^{hd-1}_M) \end{equation}
be the $hd-1$-form. This $\omega$  is the dual expression of the Jacobian matrix of some defining polynomials for the scheme
 $P(\Gamma_{\mathbb L})$. 
The crucial point of this definition is that the polynomials inside of ``d" operator,
 \begin{equation} 
\left|  \begin{array}{ccc} 
f_2(c(t_i)) & f_1(c(t_i)) & f_0(c(t_i))\\
f_2(c(t_1)) & f_1(c(t_1)) & f_0(c(t_1))\\
f_2(c(t_2)) & f_1(c(t_2)) & f_0(c(t_2))
\end{array}\right| \end{equation}
 $i=3, \cdots, hd+1$ generically define the scheme (the scheme-theorectical under $P$), 
$$P(\Gamma_{\mathbb L}).$$

Using successive blow-ups at a ``boundary point" of $P(\Gamma_{\mathbb L})$, we proved that 
\bigskip

\begin{proposition} For generic choices of quintics, $f_0, f_1, f_2$, and 
$$(t_1, \cdots, t_{hd+1})\in Sym^{hd+1}(\mathbf P^1), $$
$\omega$ is a non-zero differential form of degree $hd-1$ when restricted to a non empty open set of  $P(\Gamma_{\mathbb L}) $, i.e.
the set of global sections $\{\phi_i\}_{i=3, \cdots, hd+1}$ is linearly independent in the
$\mathcal O(P(\Gamma_{\mathbb L})  )$ module, 
$$H^0(\Omega_M\otimes \mathcal O_{P(\Gamma_{\mathbb L}) }).$$

\end{proposition}

\bigskip

Next by an argument on Zariski tangent spaces, mainly from the fact that the ideal of scheme,
$$P(\Gamma_{\mathbb L})$$ is generated by 
polynomials,
 \begin{equation} 
\left|  \begin{array}{ccc} f_2(c(t_1)) & f_1(c(t_1)) & f_0(c(t_1))\\
f_2(c(t_2)) & f_1(c(t_2)) & f_0(c(t_2))\\
f_2(c(t_i)) & f_1(c(t_i)) & f_0(c(t_i))\end{array}\right|\end{equation}
we obtain 
\bigskip

\begin{proposition}
If $\omega$ is non-zero  on $P(\Gamma_{\mathbb L})$(The algebraic equivalence of this is that the set $\{\phi_i\}_{i=3, \cdots, hd+1}$ 
is linearly independent  in the 
$\mathcal O(P(\Gamma_{\mathbb L}))$ module), then the Zariski tangent space of 
$P(\Gamma_{\mathbb L})$ at a generic maximal point must be 
\begin{equation}
dim(M)-deg(\omega)
\end{equation}
\end{proposition}

\bigskip

{\bf Remark} The form $\omega$ is not invariant under the $GL(2)$  action, but
the zero locus  $\{\omega=0\}\subset \mathbb C^{(n+1)(d+1)}$ is.\footnote{ The ideal of $\{\omega=0\}$ is a Jacobian ideal.} The form $\omega$ depends on the generic choice of $t_1, \cdots, t_{hd+1}$,
but $\{\omega=0\}$ does not.  The proof of proposition 1.2 is the main body of the section 2. It is achieved by successive blow-ups at a rational curve 
lying on  a product of $h$ distinct planes in general positions. \bigskip

Continuing from this proposition, by the surjectivity of $\Gamma_{\mathbb L}$ to $\mathbb L$, the argument on Zariski tangent 
spaces shows that the dimension of Zariski tangent space of
$\Gamma_{\mathbb L}$ at a generic point must be the same as that of Zariski tangent space of
$P(\Gamma_{\mathbb L})$ at a generic point. Thus the proposition 1.2 implies that
the dimension of the Zariski tangent space of $\Gamma_{\mathbb L}$ at a generic point is
$h+1$. 
By lemmas 2.8, 2.9,
this directly leads to 
$$H^1(N_{c_0/X_0})=0.$$
\bigskip

\bigskip

\vfill\eject

\section{Proof of Main theorem}

The proof of theorem 1.1 in length is unbalanced with other sections.  In order to be clear, we divided it into two steps that deal with two different types of
hypersurfaces. \bigskip

\subsection{Calai-Yau and  Fano hypersurfaces}\medskip\quad

The first case is $n+1-h\geq 0$. The hypersurfaces $X_0$ in this case is either 
 Fano when $n+1-h>0$ or Calabi-Yau when $n+1-h=0$. Then we can repeat the  proof in [10] to 
prove the main theorem.  However some of the steps need to be altered from [10]. In order to
insure the correctness of the proof, in the following we  go through 
the proof of [10] step by step, and make changes when they are necessary.  Also
the proof for the Fano case is identical to that for the Calabi-Yau's. So in this subsection we only prove the case when
$X_0$ is Calabi-Yau. Therefore we assume $$h=n+1.$$  We start with notations.\bigskip

\subsubsection{Technical notations}\quad\smallskip

In this section, we collect all technical notations and definitions used in section 2.1 . Some of them may already be defined before.\bigskip

{\bf Notations}:
\par
(1) $S$ denotes the space all hypersurfaces of degree $h$, i.e. $S=\mathbf P(H^0(\mathcal O_{\mathbf P^n}(h)))$.\par
Let  $[f]$ denote the image of $f$ under the map
$$\begin{array}{ccc}
H^0(\mathcal O_{\mathbf P^n}(h))-\{0\} &\rightarrow & S. \end{array}$$

(2) Let $$M$$ be
$$\mathbb C^{(n+1)(d+1)}\simeq (H^0(\mathcal O_{\mathbf P^1}(d))^{\oplus n+1}$$ 
and $M_d$ be the subset that parametrizes  regular maps $$\mathbf P^1\to \mathbf P^n$$ whose push-forward cycle has degree $d$. 
 \par
(3) Throughout the paper, if $$c:  \mathbf P^1\to \mathbf P^n, $$
is regular,   $c^\ast(\sigma)$   denotes the pull-back section  of
section $\sigma$ of some bundle over $\mathbf P^n$. The vector bundles will not always be specified, but they are apparent in the context.\par
(4) Let $Y$ be a scheme,  $y\in Y$ be a closed point, $Z\subset Y$ be a subscheme (open or closed) and $\mathcal M$ be a quasi-coherent
sheaf of $\mathcal O_{Y}$-module.
Then $\mathcal O_{y, Y}$ denotes the local ring, $\Omega_{Y}$ denotes the sheaf of differentials, 
  $\mathcal M|_{(Z)}$ denotes the pull-back sheaf module $i^\ast(\mathcal M)$ where $i: Z\hookrightarrow Y$ is the embedding. We call  
$\mathcal M|_{(Z)}$ the restriction of $\mathcal M$ to $Z$. 
$\mathcal M|_{Z}$ denotes the localization of $\mathcal M$ at $Z$, which is a $\mathcal O_{Z, Y}$ module.  Thus
$$\mathcal M|_{(\{y\})}=\mathcal M|_Z\otimes k(y),$$
where $k(y)$ is the residue field of the maximal point $\{y\}$.  
\par
If $Y$ is quasi-affine scheme, $\mathcal O (Y)$ denotes the ring of regular functions on $Y$.  

\par

(5) Let $\alpha\in T_{c_0}M$, and $$g: M\to H^0(\mathcal O_{\mathbf P^1}(r))$$
be a regular map. 
Then the  image $g_\ast(\alpha)$ of $\alpha$ under the differential map at $c_0$ is denoted by
\begin{equation}
{\partial g(c_0(t))\over \partial \alpha}\in T_{g(c_0)}(H^0(\mathcal O_{\mathbf P^1}(r)))=H^0(\mathcal O_{\mathbf P^1}(r)).
\end{equation}
(use the identification $T_{g(c_0)}(H^0(\mathcal O_{\mathbf P^1}(r)))=T_0(H^0(\mathcal O_{\mathbf P^1}(r)))$).
\par
(6) If $Y$ is a scheme, $|Y|$ denotes the induced reduced scheme of $Y$.

\bigskip

\begin{definition}   
Let
\begin{equation}
c_0^\ast(T_{X_0})\simeq \mathcal O_{\mathbf P^1}(a_1)\oplus \cdots \oplus O_{\mathbf P^1}(a_{n-1}).\end{equation}
where $$a_1\geq \cdots \geq a_{n-1},\ and\  \sum_i a_i=0$$
We'll fix an isomorphism in (2.2)  throughout. 
Let $E$ be the pull-back of the summand 
$\sum_{a_i\geq  0} \mathcal O_{\mathbf P^1}(a_i)$ under the map
$$c_0^\ast(T_{X_0})\to N_{c_0/X_0}.$$
So  $E$ is the sub-bundle of the bundle $c_0^\ast(T_{X_0})$,  generated by all the holomorphic sections of
 $c_0^\ast(T_{X_0})$.
\bigskip
\end{definition}

\bigskip

\begin{definition} 
(a) 
 If $f\in H^0(\mathcal O_{\mathbf P^n}(h))$ is a degree $d$  polynomial other than $f_0$, we denote the direction of
the line through two points $[f], [f_0]$ in the projective space, $\mathbf P(H^0(\mathcal O_{\mathbf P^n}(h))$  by $\overrightarrow f$.
So $$\overrightarrow f\in T_{[f_0]}\mathbf P(H^0(\mathcal O_{\mathbf P^n}(h)).$$\par
(b) 
Note that the vector $\overrightarrow f$ is well-defined up-to a non-zero multiple. 
In case when $c_0$ can deform to all hypersurfaces to the first order, i.e. the map in (2.7) below is surjective, this naturally gives a 
section $<\overrightarrow f>$ of the  bundle $c_0^\ast(T_{\mathbf P^n})$ (may not be unique),  to each deformation $ \overrightarrow f$ of the
hypersurface $f_0$. This is easily can be understood as the direction of the moving $c_0$ in the deformation $(\overrightarrow f, <\overrightarrow f>)$ of
the pair $(c_0, f_0)$.

\end{definition}

\bigskip

\begin{definition}
Let $\Gamma$  be an irreducible component of the incidence scheme 
\begin{equation}\{(c, f)\subset M\times \mathbf P(H^0(\mathcal O_{\mathbf P^n}(h))): c^\ast(f)=0\}
\end{equation} 
that dominates $S=\mathbf P(H^0(\mathcal O_{\mathbf P^n}(h)))$. \end{definition}
Let $(c_0, [f_0])\in \Gamma$ be a generic point.  Throughout the paper we assume that such a $\Gamma$ exists. 
\bigskip

{\bf Remark}: The existence of such a $\Gamma$ is equivalent to  the assumption of theorem 1.1: 
$X_0$ is generic. Results in section 2.1.2 only need a weaker assumption, but the main propositions 1.2, 1.3 rely on this stronger assumption---$\Gamma$ exists.  In this paper, to avoid the distraction, we use the unified and consistent assumption--$\Gamma$ exists. 
\bigskip

\begin{definition} Let $f_1, f_2 \in H^0(\mathcal O_{\mathbf P^n}(h))$ be two quintics different from $f_0$.
Let $\mathbb L$ be an open set of the plane in 
$$ \mathbf P(H^0(\mathcal O_{\mathbf P^n}(h)))$$ spanned by $[f_0], [f_1], [f_2]$ and centered around $[f_0]$. 
 \end{definition}

\bigskip

\begin{definition} Let
\begin{equation} 
\Gamma_{\mathbb L}=\Gamma\cap ( M\times \mathbb L)
\end{equation} be an irreducible component of the restriction of $\Gamma$ to $M\times \mathbb L$ such that it is onto ${\mathbb L}$, and
\begin{equation} 
\Gamma_{f_0}, \ for\ generic\ f_0\in \mathbb L
\end{equation} is an irreducible component of
$$P(\Gamma\cap ( M\times \{[f_0]\}))$$
where $P$ is the projection to $M$.

 \end{definition}

\bigskip

\begin{definition} 
Let $V$ be a smooth analytic variety with analytic coordinates $x_1, \cdots, x_p$, Let $f_1, \cdots, f_m$ be holomorphic functions on $V$.
We define

\begin{equation} \begin{array}{c} 
\left (\begin{array}{cccccc}
{\partial f_1\over \partial x_1}  &{\partial f_1\over \partial x_2}   & \cdots & {\partial f_1\over \partial x_p}   \\
{\partial f_2\over \partial x_1}  &{\partial f_2\over \partial x_2} & \cdots & {\partial f_2\over \partial x_p}  \\
 \vdots & \vdots & \cdots &\vdots \\
{\partial f_m\over \partial x_1}  &{\partial f_m\over \partial x_2} & \cdots & {\partial f_m\over \partial x_p}
\end{array}\right). \end{array}\end{equation}

to be the Jacobian matrix of functions $f_1, \cdots, f_m$. This Jacobian matrix depends on the coordinates 
$x_1, \cdots, x_p$. 

\end{definition}

\bigskip

We defined the matrix for the set of functions $f_1, \cdots, f_m$. One may wish to compare this definition with that of Jacobian ideals which
is independent choice of those functions.

\subsubsection{First order}

\smallskip\quad

 Let's start the problem in its first order.\bigskip

\begin{lemma} Let $f_0$ be a generic hypersurface containing a rational map $c_0$ as before. 
If $(c_0, [f_0])\in |\Gamma|$ is generic , then the projection
\begin{equation} \begin{array}{ccc}T_{(c_0, [f_0])}\Gamma &\stackrel{P^s}
\rightarrow & T_{[f_0]}S\end{array}\end{equation}
is surjective, where $S=\mathbf P(H^0(\mathcal O_{\mathbf P^n}(h))$.

\end{lemma}

\bigskip

\begin{proof}
Let $|\Gamma |\subset \Gamma$ be the reduced scheme of the scheme $\Gamma$.
By the genericity of $f_0$,  the projection 
\begin{equation}\begin{array} {ccc} |\Gamma | &\rightarrow &
S\end{array}\end{equation}
is dominant. Hence in a neighborhood of a generic point 
$(c_0, [f_0])\in |\Gamma|$, the projection is a smooth map. Thus 
\begin{equation}\begin{array} {ccc} T_{(c_0, [f_0])}|\Gamma | &\rightarrow &
T_{[f_0]}S\end{array}\end{equation}
is surjective. This proves the lemma

\end{proof}
\bigskip

 To elaborate definition 2.2, we apply this lemma to obtain that for any $\alpha\in T_{[f_0]}S$,  there is  
a section denoted by 
$$<\alpha>\in H^0(c_0^\ast(T_{\mathbf P^n}))$$ such that $$ (\alpha, <\alpha>)$$ is tangent to the 
universal hypersurface $$\mathcal X=\{(x, f): x\in div(f)\}$$ in $$\mathbf P^n\times S.$$
 Note that $<\alpha>$ is  unique up to a section in $H^0(c_0^\ast(T_{X_0}))$.  But we will always fix $<\alpha>$ as in definition 2.2.
 
\bigskip

\bigskip
\subsubsection{The incidence scheme}\quad\smallskip

In this subsection, we study the Zariski tangent spaces of various incidence schemes to reveal
a connection between the incidence scheme and the normal sheaf.  \bigskip

\begin{lemma} Let $[f_0]\in S$ be a generic point, $\mathbb L_1\subset S$ an open set of the pencil containing 
$f_0$ and another  quintic $f_1$.   Let $(c_0, [f_0])\in \Gamma_{\mathbb L_1}$ be generic. Then
 
\par
(a)  \begin{equation}
{T_{c_0}\Gamma_{f_0}\over ker} \simeq H^0(c_0^\ast(T_{X_0})).\end{equation}
where $ker$ is a line in $T_{c_0}\Gamma_{f_0}$. 

\par
(b) 
\begin{equation} dim( T_{(c_0, [f_0])}\Gamma_{\mathbb L_1})=dim( T_{c_0}\Gamma_{f_0})+1, 
\end{equation}
and furthermore
\begin{equation} dim( T_{c_0}P(\Gamma_{\mathbb L_1}))=dim( T_{c_0}\Gamma_{f_0})+1, 
\end{equation}
\par

\end{lemma}

\bigskip

\begin{proof}
(a). Let $a_i(c, f), i=0, \cdots, hd$ be the coefficients of 
polynomial $f(c(t))$ in parameter $t$. Then the scheme
$$\Gamma$$ in $M\times \mathbf P^n$ is defined by homogeneous
polynomials $$a_i(c, f)=0, i=0, \cdots, hd, \ locally.$$ 
 Let $\alpha\in T_{c_0}M$. The equations
\begin{equation} {\partial  a_i(c_0, f_0)\over \partial \alpha}=0, \ all\ i \end{equation}  by the definition,  are necessary and sufficient conditions
for $\alpha$ lying in $$T_{(c_0, [f_0])}\Gamma_{f_0}.$$ 
On the other hand there is an evaluation map $e$:
\begin{equation}\begin{array}{ccc}
M\times \mathbf P^1 &\rightarrow & \mathbf P^n\\
(c, t) &\rightarrow & c(t)
\end{array}\end{equation}
The differential map
$e_\ast$ gives a morphism $e_m$:
\begin{equation}\begin{array}{ccc}
T_{c_0}M &\stackrel{e_m}\rightarrow & H^0(c_0^\ast( T_{\mathbf P^n}))\\
\alpha &\rightarrow & e_\ast (\alpha)
\end{array}\end{equation}

Suppose there is an $\alpha$ such that
$e_\ast (\alpha)=0$. We may assume $c_0$ is a map
$$\mathbb C^1\to \mathbb C^{n+1}-\{0\}.$$
 
Since $c_0$ is birational to its image, there is a Zariski open set  $$U_{\mathbf P^1}\subset \mathbf P^1$$
and an open set $$V\subset \{(c_0(t))\}\subset \mathbb C^{n+1}-\{0\}$$
such that $c_0'|_{U_{\mathbf P^1}}$ is an isomorphism 
$$U_{\mathbf P^1}\to V, $$ where
$c_0'$ is the map from $U_{\mathbf P^1}$ to $\mathbb C^{n+1}$ induced from $c_0$.
Due to the equation $e_\ast (\alpha)=0$, on $T_{t}U_{\mathbf P^1}$
$$(\alpha_0(t), \cdots, \alpha_n(t))=\lambda(t) c_0(t)$$
on $V$ (at each point $(c_0(t), \cdots, c_n(t))$ of $V$)
where $\lambda(t)$ lies in $\mathcal O(U_{\mathbf P^1})$. Because
$(\alpha_0(t), \cdots, \alpha_n(t))$ is parallel to
$c_0(t)$ at all points $t\in \mathbf P^1$, $\lambda$ can be extended to $\mathbf P^1$. Hence
$\lambda(t)$ is in $\mathcal O(\mathbf P^1)$.  So it is a constant (independent of $t$).
Therefore $\alpha\in \mathbb C^{h(d+1)}$ is parallel to $$c_0\neq 0\in  \mathbb C^{h(d+1)}.$$
This shows that 
$$dim(ker(e_m))=1.$$
(this does not hold if $c_0$ is a multiple cover map). By the dimension count, $e_m$ must be surjective.
 For any $\alpha\in c_0^\ast( T_{\mathbf P^n})$, 
$\alpha\in c_0^\ast(T_{X_0})$ if and only if
\begin{equation} {\partial  f_0(c_0(t))\over \partial \alpha}=0, \end{equation} 
for generic $t\in \mathbf P^1$. Notice equations (2.13) and (2.16) are exactly the same. Therefore
$e_m$ induces an isomorphism 

\begin{equation} \begin{array}{ccc} {T_{c_0}\Gamma_{f_0}\over ker(e_m)} &\stackrel{e_m}\rightarrow &  H^0(c_0^\ast(T_{X_0}))\\
\end{array}\end{equation}
This proves part (a).\par

(b). Let $t_1, \cdots, t_{hd+1}\in \mathbf P^1$ be $hd+1$ distinct points of $\mathbf P^1$.

Notice  $c^\ast(f)=0$ if and only if $c^\ast(f)|_{t_i}=0$ for all $i$. 
Then $$\Gamma_{\mathbb L_1}$$ is an irreducible component of 
\begin{equation} \{ (c, f)\in M\times \mathbb L_1: c^\ast(f)|_{t_i}=0, i=1, \cdots, hd+1\}. \end{equation}
surjective to $S$ in first order. 
Let $\alpha\in T_{c_0}M$ and $\beta=\overrightarrow f_1$. Let  $\epsilon\in \mathbb C$, then
 $\epsilon \beta\in T_{[f_0]}\mathbb L_1$. 
Then the Zariski tangent space $T_{c_0}\Gamma_{\mathbb L_1}$ is
$$\{(\alpha, \epsilon \beta)\in T_{(c_0, [f_0])}(M\times \mathbb L_1):
\epsilon f_1(c_0(t_i))+{\partial f_0(c_0(t_i))\over \partial \alpha} =0, i=1, \cdots, hd+1\}$$
Because $[f_0]$ is a generic point of $\mathbb L_1$, the map (2.7) is surjective.  Thus
there is $\alpha_0\in T_{c_0}(M)$ such that 
$$ f_1(c_0(t))+{\partial f_0(c_0(t))\over \partial \alpha_0}=0$$ for all $t\in\mathbf P^1$. 
Thus $T_{(c_0, [f_0])}\Gamma_{\mathbb L_1}$ is isomorphic to
$$\{(\alpha, \epsilon)\in\mathbb C^{(n+1)(d+1)+1}:
{\partial f_0(c_0(t_i))\over \partial (\alpha-\epsilon\alpha_0)}=0, i=1,\cdots, hd+1\}.$$

Notice that the subspace with $\epsilon=0$, 
$$\{\alpha\in T_{c_0}M:
{\partial f_0(c_0(t_i))\over \partial \alpha} =0, i=1,\cdots, hd+1\},$$
is  the  tangent space of the scheme
$$\Gamma_{f_0}.$$
Thus the dimensions of them differ by $1$. 

\par

Next we prove the assertion for the scheme-theoretical image $P(\Gamma_{\mathbb L_1})$. 
Using  an expression of $\Gamma_{\mathbb L_1}$, $P(\Gamma_{\mathbb L_1})$ is defined
by polynomials
\begin{equation}
f_1(c(t_i))f_0(c(t_j))-f_1(c(t_j))f_0(c(t_i))=0
, 1\leq i, j\leq hd+1.\end{equation}
Assume non of $t_i, i=1, \cdots, hd+1$ is a zero of
$f_1(c_0(t))=0$. Then there exists an open set $U$ of $M$ around $c_0$ such that
\begin{equation} P(\Gamma_{\mathbb L_1})\cap U\end{equation}
is defined by
$hd$ equations
\begin{equation}
f_1(c(t_{hd+1}))f_0(c(t_j))-f_1(c(t_j))f_0(c(t_{hd+1}))=0, j=1, \cdots, hd.\end{equation}
We'll denote $U$ by $M$.
Then the Zariski tangent space $T_{c_0} P(\Gamma_{\mathbb L_1})$ is defined by
\begin{equation}\begin{array} {cc} &
  f_1(c_0(t_{hd+1})){\partial f_0(c_0(t_j))\over \partial \alpha}-f_1(c_0(t_j)){\partial f_0(c_0(t_{hd+1}))\over \partial \alpha}=0, \\
& j=1, \cdots, hd.\end{array}
\end{equation}

where $\alpha\in T_{c_0}M$.
This is the same as
\begin{equation}
{\partial f_0(c_0(t_j))\over \partial \alpha}-{f_1(c_0(t_j))\over f_1(c_0(t_{hd+1})) }{\partial f_0(c_0(t_{hd+1}))\over \partial \alpha}=0
\end{equation}

Next we view ${\partial f_0(c_0(t_j))\over \partial \alpha}$ as an element in
$$(T_{c_0}M)^\ast. $$
If $${\partial f_0(c_0(t_i))\over \partial \alpha}, i=1, \cdots, hd+1$$ are linearly independent, then
$dim(T_{c_0}\Gamma_{f_0})=(n+1-h)d+n$ and by (2.23), $$dim(T_{c_0}P(\Gamma_{\mathbb L_1}))=(n+1-h)d+n+1. $$ The lemma is proved.

If $${\partial f_0(c_0(t_i))\over \partial \alpha}, i=1, \cdots, hd+1$$ are linearly dependent, there are two cases:\par
(1) All solutions $\alpha_0$  to (2.23) satisfy $${\partial f_0(c_0(t_{hd+1}))\over \partial \alpha}=0.$$
(2) Some solutions $\alpha_0$ to (2.23) do not satisfy $${\partial f_0(c_0(t_{hd+1}))\over \partial \alpha}=0.$$
The case (1) is false. Because if $${\partial f_0(c_0(t_{hd+1}))\over \partial \alpha_0}=0,$$
for all solution $\alpha_0$ of $(2.23)$, 
then ${\partial f_0(c_0(t_j))\over \partial \alpha_0}=0$ for all $j=1, \cdots, hd+1$. Hence all
solutions $\alpha_0$ to (2.23) must be  sections of $E$ (see definition 2.1 for $E$). As we know 
$$<\overrightarrow f_1>$$ is a solution to (2.23), but 
$${\partial f_0(c_0(t))\over \partial <\overrightarrow f_1>}=-f_1(c_0(t))$$ 
which means $<\overrightarrow f_1>$ is not a section of $E$. This is a contradiction. 

\par
In case (2), the solutions $\alpha_0$ to (2.23) must either satisfy 
$${\partial f_0(c_0(t_j))\over \partial \alpha}=0$$ for all $j=1, \cdots, hd+1$, in which case, they are sections 
of $E$, or are uniquely expressed as
\begin{equation}
{\partial f_0(c_0(t_j))\over \partial \alpha_0}={f_1(c_0(t_j))\over f_1(c_0(t_{hd+1})) }{\partial f_0(c_0(t_{hd+1}))\over \partial \alpha_0}, 
j=1, \cdots, hd+1
\end{equation}
which are exactly $<\overrightarrow f_1>$. 
Thus the vector $<\overrightarrow f_1>$ offers another dimension to 
$$dim(T_{c_0}P(\Gamma_{\mathbb L_1})).$$
So 
\begin{equation}
T_{c_0}P(\Gamma_{\mathbb L_1})=T_{c_0}\Gamma_{f_0}\oplus \mathbb C<\overrightarrow f_1>_M,
\end{equation}
where $<\overrightarrow f_1>_M\in e_m^{-1} (<\overrightarrow f_1>)$. 
The lemma is proved.
\bigskip

\end{proof}

\bigskip

Now we can describe the case when this is based on 2 dimensional plane in $S$. 
Recall $\mathbb L$ is an open set of a plane spanned by $f_0, f_1, f_2$, and $$\mathbb L_1\subset \mathbb L$$ is a pencil 
containing $[f_0]$.  
\bigskip

\begin{lemma}
For generic $c_g\in P(\Gamma_{\mathbb L_1})\subset P(\Gamma_{\mathbb L})$,
\begin{equation} 
dim( T_{c_g}P(\Gamma_{\mathbb L}))=dim( T_{c_g}P(\Gamma_{\mathbb L_1}))+1
\end{equation}

\end{lemma}

\bigskip

\begin{proof}  Consider an open set $U_{P(\Gamma_{\mathbb L})}$ of ${P(\Gamma_{\mathbb L})}$ centered at $c_g$ such that
vectors in $\mathbb C^3$, 
$$\begin{array}{cc} & \biggl( f_2(c(t_1)), f_1(c(t_1)),  f_0(c(t_1))\biggr)\\&
\biggl(f_2(c(t_2)),  f_1(c(t_2)), f_0(c(t_2))\biggr)\end{array}$$
and vectors in $\mathbb C^{hd+1}$
\begin{equation} 
\left (\begin{array}{cc}  &f_2(c(t_{hd+1}) \\&
\vdots \\ &
f_2(c(t_2)) \\
 & f_2(c(t_1)) \end{array}\right),  \left (\begin{array}{cc}  &f_1(c(t_{hd+1}) \\&
\vdots \\ &
f_1(c(t_2)) \\
 & f_1(c(t_1)) \end{array}\right)
\end{equation}
are linearly independent for all $c\in U_{P(\Gamma_{\mathbb L})}$. 
If $$\left|  \begin{array}{ccc} f_2(c(t_1)) & f_1(c(t_1)) & f_0(c(t_1))\\
f_2(c(t_2)) & f_1(c(t_2)) & f_0(c(t_2))\\
f_2(c(t_i)) & f_1(c(t_i)) & f_0(c(t_i))\end{array}\right|=0,$$

$$\biggl(f_2(c(t_i)),  f_1(c(t_i)), f_0(c(t_i))\biggr)$$ for all $i=1, \cdots, hd+1$ are linear combinations of 
$$\begin{array}{cc} & \biggl( f_2(c(t_1)), f_1(c(t_1)),  f_0(c(t_1))\biggr)\\& 
\biggl(f_2(c(t_2)),  f_1(c(t_2)), f_0(c(t_2))\biggr).\end{array}$$
Hence 
$$\left|  \begin{array}{ccc} f_2(c(t_i)) & f_1(c(t_i)) & f_0(c(t_i))\\
f_2(c(t_j)) & f_1(c(t_j)) & f_0(c(t_j))\\
f_2(c(t_l)) & f_1(c(t_l)) & f_0(c(t_l))\end{array}\right|=0,$$
for all $i, j, l$ between $1$ and $hd+1$. 

Since $$\left|  \begin{array}{ccc} f_2(c(t_i)) & f_1(c(t_i)) & f_0(c(t_i))\\
f_2(c(t_j)) & f_1(c(t_j)) & f_0(c(t_j))\\
f_2(c(t_l)) & f_1(c(t_l)) & f_0(c(t_l))\end{array}\right|=0,$$
for all $i, j, l$ between $1$ and $hd+1$ define $U_{P(\Gamma_{\mathbb L})}$, 
$hd-1$ equations
$$\left|  \begin{array}{ccc} f_2(c(t_1)) & f_1(c(t_1)) & f_0(c(t_1))\\
f_2(c(t_2)) & f_1(c(t_2)) & f_0(c(t_2))\\
f_2(c(t_i)) & f_1(c(t_i)) & f_0(c(t_i))\end{array}\right|=0,$$
for $i=3, \cdots, hd+1$
define the scheme $U_{P(\Gamma_{\mathbb L})}$. 
Then by the definition there is $f_g\in \mathbb L$ such that
\begin{equation} 
c_g^\ast(f_g)=0. 
\end{equation}

We denote $U_{\mathbb L}$ by  $\mathbb L$.  Then
 $P(\Gamma_{\mathbb L}))$ is defined by the polynomial equations

\begin{equation} 
\left|  \begin{array}{ccc} f_2(c(t_i)) & f_1(c(t_i)) & f_0(c(t_i))\\
f_2(c(t_1)) & f_1(c(t_1)) & f_0(c(t_1))\\
f_2(c(t_2)) & f_1(c(t_2)) & f_0(c(t_2))\end{array}\right| =0\end{equation}
 $i=3, \cdots, hd+1$.
\par

We may assume that $c_g=c_0$ lies in $f_0$ (by choosing appropriate basis $f_0, f_1, f_2$ of $\mathbb L$).  By  (2.29),
the Zariski tangent space $T_{c_0}P(\Gamma_{\mathbb L})$ is defined by equations
\begin{equation} 
\left|  \begin{array}{ccc} f_2(c_0(t_i)) & f_1(c_0(t_i)) & {\partial f_0(c_0(t_i))\over\partial \alpha}\\
f_2(c_0(t_1)) & f_1(c_0(t_1)) & {\partial f_0(c_0(t_1))\over\partial \alpha}\\
f_2(c_0(t_2)) & f_1(c_0(t_2)) & {\partial f_0(c_0(t_2))\over\partial \alpha}\end{array}\right| =0\end{equation}
 where $i=3, \cdots, hd+1$, $\alpha\in T_{c_0}P(\Gamma_{\mathbb L})$.
This is equivalent to that 
the column vectors of 
\begin{equation} 
\left (\begin{array}{ccc} f_2(c_0(t_{hd+1}) & f_1(c_0(t_{hd+1})) & {\partial f_0(c_0(t_{hd+1}))\over\partial \alpha_0}\\
\vdots & \vdots &\vdots\\
f_2(c_0(t_2)) & f_1(c_0(t_2)) & {\partial f_0(c_0(t_{2}))\over\partial \alpha_0}\\
f_2(c_0(t_1)) & f_1(c_0(t_1)) & {\partial f_0(c_0(t_{1}))\over\partial \alpha_0}\end{array}\right) \end{equation}
 are linearly dependent,  for some  $\alpha_0\in T_{c_0}P(\Gamma_{\mathbb L})$. 
Since 
\begin{equation} 
\left (\begin{array}{cc}  &f_2(c_0(t_{hd+1}) \\&
\vdots \\ &
f_2(c_0(t_2)) \\
 & f_2(c_0(t_1)) \end{array}\right),  \left (\begin{array}{cc}  &f_1(c_0(t_{hd+1}) \\&
\vdots \\ &
f_1(c_0(t_2)) \\
 & f_1(c_0(t_1)) \end{array}\right)
\end{equation}
are linearly independent,
there are complex numbers $\epsilon_i, i=1, 2$ such that
\begin{equation} 
\left (  \begin{array}{c}  {\partial f_0(c_0(t_{hd+1}))\over\partial \alpha_0}\\
 \vdots\\
  {\partial f_0(c_0(t_{2}))\over\partial \alpha_0}\\
 {\partial f_0(c_0(t_{1}))\over\partial \alpha_0}\end{array}\right)=\epsilon_2 \left (  \begin{array}{c}  f_2(c_0(t_{hd+1})\\
 \vdots\\
 f_2(c_0(t_2))\\
 f_2(c_0(t_1))\end{array}\right)+\epsilon_1 \left (  \begin{array}{c} f_1(c_0(t_{hd+1})\\
 \vdots\\
f_1(c_0(t_2))\\
 f_1(c_0(t_1))\end{array}\right)
\end{equation}

Let $f_3=\sum_{i=1}^2 \epsilon_i f_i$ where $\epsilon_i$ are fixed.   We may assume that $\mathbb L_1$ is the pencil containing $f_0, f_3$.  The equation (2.33) becomes
$$ {\partial f_0(c_0(t))\over \partial \alpha_0}-f_3(c_0(t))=0, for\ all\ t\in\mathbf P^1.$$

Then just as in the (2.25), 
\begin{equation}
T_{c_0}\Gamma_{\mathbb L}\simeq T_{c_0}\Gamma_{f_0}\oplus \mathbb C<\overrightarrow{f_1}>_M\oplus \mathbb C<\overrightarrow{f_2}>_M.
\end{equation}
and \begin{equation}
T_{c_0}\Gamma_{\mathbb L_1}\simeq T_{c_0}\Gamma_{f_0}\oplus \mathbb C<\overrightarrow{f_3}>_M.
\end{equation}

Then the lemma follows.

\end{proof}

\subsubsection{Space of rational curves, $M$}

\smallskip\quad

The main purpose of this section is to introduce analytic coordinates of a neighborhood of $M$, which will be used as local coordinates 
of  the blow-up $\tilde M$ of $M$ in the computation of a Jacobian matrix (2.87).  These coordinates identify a neighborhood of $\tilde M$ outside of exceptional
divisor with a neigborhood of $M$.  They are crucial.They are crucial.

\bigskip

Let $\tilde c_2=(\tilde c_2^0(t), \cdots, \tilde c_2^n)\in M$ with $$\tilde c_2^i\neq 0\in H^0(\mathcal O_{\mathbf P^1}(d)), i=0, \cdots, n.$$
We assume $\tilde c_2^i(t)=0, i\leq n-2, $ have  $(n-1)d$ distinct  zeros $$\tilde \theta_i^j, for\ i\leq n-2.$$
Then the first $n-1$ components, $ H^0(\mathcal O_{\mathbf P^1}(d))$ of 
$$ M=H^0(\mathcal O_{\mathbf P^1}(d))^{\oplus n+1}$$
 has local analytic coordinates \begin{equation}
r_l, \theta_i^j, j=1, \cdots, d, l=0, \cdots, n-2, i=0, \cdots, n-2\end{equation}
with $ \ r_i\neq 0$  around
$\tilde c_2^i$ such that
\begin{equation} c^i(t)=r_i \prod _{j=1}^d (t-\theta_i^j),  
\end{equation}
and last two components, $ H^0(\mathcal O_{\mathbf P^1}(d))$ of 
$$ M=H^0(\mathcal O_{\mathbf P^1}(d))^{\oplus n+1}$$ 
have affine coordinates $r_i, \theta_i^j, i=n-1, n$  near $\tilde c_2$, 
\begin{equation}  r_{n-1}=c_{n-1}^0, r_{n}=c_n^0, \theta_i^j={c_i^j\over c_i^0}, for \ i=n-1, n, j=1, \cdots, d
\end{equation}
where $c_i^j$ are coefficients of $c^i(t)$. Assume the values of coordinates $$c_i^0, i=n-1, n$$ for
the point $\tilde c_2$ are non-zeros, i.e. $\tilde c_2$ lies in the neighborhood of this coordinate system:
$$r_l, \theta_i^j, \ for  \ 0 \leq l, i\leq  n, j=1, \cdots, d, r_l\neq 0.$$
Let a part of coordinates values for $\tilde c_2$ be
$$\begin{array}{c} r_l=y_l, \theta_i^j=\tilde\theta_i^j, l=0, \cdots, n-2, i=0, \cdots, n, j=1, \cdots, d\\
c_{n-1}^0=y_{n-1}, c_n^0=y_n.\end{array}$$
Let $q$ be a generic, homogeneous quadratic polynomial in $z_0, \cdots, z_n$.
Let 
\begin{equation} 
h(c, t)=\delta_1q(c(t))+\delta_2 c_{n-1}(t) c_n(t).
\end{equation}
for $c\in M$, where $\delta_i, i=1, 2$ are two none zero complex numbers.
   Let $\beta_1, \cdots, \beta_{2d}$ be the zeros of 
$h(\tilde c_2, t)=0$.
  Furthermore we  assume $\beta_i, i=1, \cdots, 2d$ are distinct and non-zeros for generic $\delta_i$ and specific
$q=z_1\cdots z_{n-2}$.

\bigskip

\begin{proposition}
Let $U_{\tilde c_2}\subset M$ be an analytic neighborhood of $\tilde c_2$. 

Let \begin{equation}\begin{array}{ccc}
g: U_{\tilde c_2} &\rightarrow &  \mathbb C^{(n+1)(d+1)}
\end{array}\end{equation}
be a regular map that is defined
by 
\begin{equation}\begin{array}{cc} &
g (\theta_0^1, \cdots, \theta_n^d, r_0, \cdots, r_n) \\
 &\|\\ &
(\theta_0^1, \cdots, \theta_{n-2}^d, h(c, \beta_1), \cdots, h(c, \beta_{2d}), 
r_0, \cdots, r_n).
\end{array}\end{equation} 

Then $g$ is an isomorphism to its image.

\end{proposition}
\bigskip

\begin{proof}
It suffices to prove the differential of $g$ at $\tilde c_2$ is an isomorphism for a SPECIFIC $q$. So we assume  that 
$$ \delta_1=\delta_2=1, q=z_1\cdots z_{n-2}$$
This is a straightforward calculation of the Jacobian of $g$. We may still assume that $\beta_i, i=1, \cdots, 2d$ are distinct.
The Jacobian 
\begin{equation}
{\partial g(\tilde \theta_0^1, \cdots, \tilde \theta_{n-2}^d, y_0,\cdots, y_n,  h(\tilde c_2, \beta_1),\cdots, h(\tilde c_2, \beta_{2d})
\over \partial (\theta_0^1, \cdots, \theta_{n-2}^d, r_0, \cdots, r_n, \theta_{n-1}^1, \cdots, \theta_n^d)}
\end{equation} is equal to 
another Jacobian
\begin{equation}
Ja=\left|\begin{array}{ccc}
{\partial  h(c, \beta_1)\over \partial \theta_{n-1}^1}  &\cdots  & {\partial h(c, \beta_{1})\over \partial\theta_n^d }\\
\vdots & \vdots & \vdots\\
{\partial  h(c, \beta_{2d})\over \partial \theta_{n-1}^1}  &\cdots  & {\partial h(c, \beta_{2d})\over \partial\theta_n^d }
\end{array}\right|_{\tilde c_2}
\end{equation}
Let $T_i, i=0, d$   be the determinant
\begin{equation}
\left|\begin{array}{ccc}\beta_{i+1}  & \cdots & \beta_{i+1}^d\\

\vdots & \vdots & \vdots \\
\beta_{i+d}  & \cdots & \beta_{i+d}^d
\end{array}\right|
\end{equation}
Then we compute the determinant to have 
\begin{equation}
Ja=(-1)^{d}\delta_2  T_0T_d \prod_{i=1}^d (\tilde c_2^{n-1}(\beta_{d+i})\tilde c_2^n(\beta_{i})-  \tilde c_2^{n-1}(\beta_{i})\tilde c_2^n(\beta_{d+i})).\end{equation}
Since $\beta_i$ are distinct and non-zeros, $$T_0\neq 0, T_d\neq 0.$$

Since ${\tilde c_2^{n-1}(t)\over \tilde c_2^n(t)}$ is a rational function and
$$deg(\tilde c_2^{n-1}(t))=deg(\tilde c_2^n(t))=d$$
then we can always arrange the subscripts $i$  of $\beta_i$ such that each number 
\begin{equation} ({\tilde c_2^{n-1}(\beta_{d+i})\over \tilde c_2^n(\beta_{d+i})}-{\tilde c_2^{n-1}(\beta_{i})\over \tilde c_2^n(\beta_{i})})
\end{equation} is not zero. Hence
$$ \prod_{i=1}^d (\tilde c_2^{n-1}(\beta_{d+i})\tilde c_2^n(\beta_{i})-  \tilde c_2^{n-1}(\beta_{i})\tilde c_2^n(\beta_{d+i}))\neq 0$$
Thus $Ja$ is non-zero.

Therefore 
\begin{equation}
{\partial g(\tilde \theta_0^1, \cdots, \tilde \theta_{n-2}^d, y_0,\cdots, y_n,  h(\tilde c_2, \beta_1),\cdots, h(\tilde c_2,  \beta_{2d}))
\over \partial (\theta_0^1, \cdots, \theta_{n-2}^d, r_0, \cdots, r_n, \theta_{n-1}^1, \cdots, \theta_n^d)}\neq 0
\end{equation}

We complete the proof. 

\end{proof}

\bigskip

\begin{definition}
Let $\epsilon_i=h(c, \beta_i), i=1, \cdots, 2d$. Then by proposition 2.10, 
\begin{equation} \theta_0^1, \cdots, \theta_{n-2}^d, r_0, \cdots, r_n, \epsilon_1, \cdots, \epsilon_{2d}
\end{equation}
are the local analytic coordinates of $M$ around $\tilde c_2$, and
$\tilde c_2$ corresponds to the coordinate values
\begin{equation}\begin{array} {c}
 \theta_i^j=\tilde \theta_i^j, i\leq n-2 , j=1, \cdots, d\\
r_l=y_l\neq 0, l=0, \cdots, n\\
\epsilon_i=0, i=1, \cdots, 2d
\end{array}\end{equation}

\end{definition}

\bigskip

\subsubsection{Differential sheaf}

\smallskip\quad

Next we prove theorem 1.1, i.e. \bigskip

\begin{equation} H^1(N_{c_0/X_0})=0\end{equation}
at generic $(c_0, [f_0])\in \Gamma$ for Calabi-Yau $X_0$.

\bigskip

\bigskip

Choose a  homogeneous coordinate system $[z_0, \cdots, z_n]$ for $\mathbf P^n$. 
Let 
\begin{equation} 
f_3=z_0\cdots z_{n-2}(\delta_1q+ \delta_2z_{n-1} z_n).
\end{equation}
where $\delta_i$ are two complex non-zero complex numbers, and $q$ is a generic quadric in $\mathbf P^n$. 
Let 
$\tilde c_2\in M$ and 
$$f_3(\tilde c_2( t))\neq 0.$$
As before we denote the zeros of $\tilde c_2^i(t)=0$ by $\tilde\theta_i^j, j\leq n-2$ and zeros of
\begin{equation}
(\delta_1q+ \delta_2z_{n-1} z_n|_{\tilde c_2(t)})=0
\end{equation}
by $\beta_i, i=1, \cdots, 2d$. 
We assume $\tilde \theta_i^j, j\leq n-2$ are distinct, and $\beta_i, i=1, \cdots, 2d$ are also distinct. 

\bigskip

\begin{lemma} 
Let $t_1, \cdots, t_{hd}$ be the zeros of  $f_3(\tilde c_2( t))$.   Recall 
in definition 2.11, 
$$\begin{array} {c}
 \theta_i^j=\tilde \theta_i^j, i\leq n-2, j=1, \cdots, d\\
r_l=y_l\neq 0, l=0, \cdots, n\\
\epsilon_i=0, i=1, \cdots, 2d
\end{array}$$ are analytic coordinates of $M$ around the point $\tilde c_2$.

Then the Jacobian matrix

\begin{equation}\begin{array} {c} J( \tilde c_2)\\
\|\\
\left(\begin{array}{cccccc}
{\partial f_3(\tilde c_2(t_{1}))\over \partial \theta_0^1}  &\cdots &
{\partial f_3(\tilde c_2(t_{1}))\over \partial \theta_{n-2}^d} & {\partial f_3(\tilde c_2(t_{1}))\over \partial \epsilon_1} &\cdots &
{\partial f_3(\tilde c_2(t_{1}))\over \partial \epsilon_{2d}}\\
{\partial f_3(\tilde c_2(t_{2}))\over \partial \theta_0^1}  &\cdots &
{\partial f_3(\tilde c_2(t_{2}))\over \partial \theta_{n-2}^d} & {\partial f_3(\tilde c_2(t_{2}))\over \partial \epsilon_1} &\cdots &
{\partial f_3(\tilde c_2(t_{2}))\over \partial \epsilon_{2d}}\\
{\partial f_3(\tilde c_2(t_{3}))\over \partial \theta_0^1}  &\cdots &
{\partial f_3(\tilde c_2(t_{3}))\over \partial \theta_{n-2}^d} & {\partial f_3(\tilde c_2(t_{3}))\over \partial \epsilon_1} &\cdots &
{\partial f_3(\tilde c_2(t_{3}))\over \partial \epsilon_{2d}}\\
\vdots&\vdots &\vdots &\vdots &\vdots&\vdots\\
\vdots&\vdots &\vdots &\vdots &\vdots&\vdots\\
{\partial f_3(\tilde c_2(t_{hd}))\over \partial \theta_0^1}  &\cdots &
{\partial f_3(\tilde c_2(t_{hd}))\over \partial \theta_{n-2}^d} & {\partial f_3(\tilde c_2(t_{hd}))\over \partial \epsilon_1} &\cdots &
{\partial f_3(\tilde c_2(t_{hd}))\over \partial \epsilon_{2d}}
\end{array}\right)\end{array}\end{equation}
is equal to 
a diagonal matrix $D$ whose diagonal entries are
\begin{equation}\begin{array}{c}
{\partial f_3(\tilde c_2(t_1))\over \partial \theta_0^1}, \cdots, {\partial f_3(\tilde c_2(t_{(n-1)d}))\over \partial \theta_{n-2}^d}, \\
{\partial f_3(\tilde c_2( t_{(n-1)d+1}))\over \partial \epsilon_1}, \cdots, 
{\partial f_3(\tilde c_2( t_{hd}))\over \partial \epsilon_{2d}}.\end{array}
\end{equation}
which are all non-zeros.
\end{lemma}

\bigskip

\begin{proof}  Note $\tilde \theta_i^j$ are distinct and $\beta_i$ are also distinct by the genericity of $q$.  Thus the coordinates
in definition 2.11 exist. 
It suffices to show all non diagonal entries of (2.53) are zeros. We can rewrite 
\begin{equation}
f_3(c(t))=y \prod_{j=1}^{hd}(t-\alpha_j)
\end{equation}
around $\tilde c_2$.
Then $y, \alpha_j$ are all functions of the analytic coordinates in definition 11, 
\begin{equation}
\theta_i^j, \epsilon_1, \cdots, \epsilon_{2d}, r_l.
\end{equation}
(Next we simply show that $\alpha_i$ correspond to coordinates $\theta_i^j, \epsilon_n, i\leq n-2$).
By the definition the Jacobian matrix of  
\begin{equation}
{\partial (\alpha_1, \cdots,  \alpha_{hd})\over \partial (\theta_0^1, \cdots, \theta_{n-2}^d, \epsilon_1,\cdots,\epsilon_{2d})}
\end{equation}
is a  non-zero diagonal matrix when evaluated at $\tilde c_2$.  Since it is clear that
$$
{\partial f_3 (\tilde c_2(t_i))\over \partial \alpha_j}=0, i\neq j, 
$$
the non-diagonals of $J(\tilde c_2)$ are also $0$. 
Upto a non-zero constant, the diagonal entries $J(\tilde c_2)$ are 
$$\prod_{j\neq 1}(t_1-t_j), \prod_{j\neq 2}(t_2-t_j), \cdots, \prod_{j\neq hd}(t_{hd}-t_j).$$
We complete the proof. 
\end{proof}

\bigskip

\subsubsection{\bf Non-vanishing $hd-1$-form $\omega$ }
\quad\smallskip

The section 2.1.4 and lemma 2.12 are just the preparation for the proof. Section 2.1.3 is the one that is meant to dig into the problem, but it is
only to the first order. With only the results in section 2.1.3, we can't go far because theorem 1.1 touches upon the higher orders of deformations of pairs.
The following lemma is the reflection of this philosophy.\bigskip

\begin{lemma} 

The $hd$-$1$ form $\omega$ defined in (1.5) is a non-zero form when it is evaluated at generic points of $P(\Gamma_{\Bbb L})$, i.e.
the reduction $\bar \omega$ in the module, 
$$H^0(\Omega_M\otimes \mathcal O_{P(\Gamma_{\Bbb L})})$$ is non zero.

\end{lemma}
\bigskip

This is the proposition 1.2.

\bigskip

It suffices to prove lemma 2.13 for a special choice of $f_0, f_1, f_2$, and we only need to produce one point (any point) on $P(\Gamma_{\mathbb L})$, 
at which $\omega$ is non-vanishing.\par
So let $z_0, z_1, \cdots, z_n$ be general homogeneous coordinates of
$\mathbf P^n$. Let $$f_2=z_0\cdots z_n.$$  
Let 
$$f_1=z_0\cdots z_{n-2} q,$$
where $q$ is a generic quadratic polynomial in $z_0, \cdots, z_n$. 
   Choose another generic $f_0$.   We obtain an open set $\mathbb L_1$ of pencil through $f_0, f_2$, and an
 open set $\mathbb L$ of 2-dimensional plane containing $f_0, f_1, f_2$ in $S$.
 Let $P(\Gamma_{\mathbb L})$ and $P(\Gamma_{\mathbb L_1})$ be as defined in lemmas 2.8, 2.9. We choose $P(\Gamma_{\mathbb L_1})$ to be irreducible, and to 
be contained in $P(\Gamma_{\mathbb L})$ for generic $f_0$. 
We may assume a generic point $c=(c^0, \cdots, c^n)\in  P(\Gamma_{\mathbb L_1})$ does not have multiple zeros with coordinates planes, i.e.
$c^i=0, i=0, \cdots$ have $hd$ distinct roots (This is because we can always choose a generic coordinates' system $z_i$ with respect to the fixed $c_0$). 
Let $$c_2\in P(\Gamma_{\mathbb L_1})\subset P(\Gamma_{\mathbb L})$$ such
that $f_2(c_2(t))=0$.   By the genericity of $f_0$ and $z_0$, 
we assume 
$$P(\Gamma_{\mathbb L_1})\not\subset \{z_0=0\}.$$
Thus we may assume $c_2=[0, c^1_2, \cdots, c_2^n]$ where
$c_2^i, i\neq 0$ are non-zero sections of $H^0(\mathcal O_{\mathbf P^1}(d))$. 
Since the proofs for other cases with more zero sections $c_2^i$ are the same,  it suffices for us to prove this case only. 
It is not difficult to see $\omega$ is zero at $c_2$ for the choice
of $\mathbb L$.    But we would like to show that $\omega$ is  not identically zero on  $P(\Gamma_{\mathbb L})$.    The technique is ``blow-up".  \par
Let's  first describe it in the most general term: we would
construct a birational map, the composition of successive blow-ups, 
\begin{equation}\begin{array}{ccc}
\pi: Y &\rightarrow &P(\Gamma_{\mathbb L})
\end{array}\end{equation}
which is an isomorphism on
$Y-\pi^{-1}(B')$ where $B'$ is a proper closed sub-scheme of $P(\Gamma_{\mathbb L})$ that contains $c_2$, and
$\pi^{-1}(B')$ is also a proper closed sub-scheme.
We then compute to find out that 
\begin{equation}
\pi^\ast(\omega)= g \omega'
\end{equation}
on $Y$ where $g$ is a non-zero rational function on $Y$, and  $\omega'$ does not vanish at a point $p\in \pi^{-1}(B')$. Because $\pi$ is birational,
$\omega$ is a non-zero form restricted to  $P(\Gamma_{\mathbb L})$. The key to this assertion is that
$\omega'$ is non-zero at one point $p$. 
The process of blow-ups is messy and lengthy. In order to organize them,
we  divide them into two different types: \par
(1) blow-ups used to resolve the vanishing of sections $c_2^i, i=0, \cdots, n$ in the rational map $c_2$. These will be only used
in step 1 below.\par
(2) blow-ups used to resolve the multiple zeros of $\tilde c_2$ with coordinates planes. This
is the case where $c_2$ could be a constant map. These will be  used in step 2 below. 

\par

Let's see the details.  
Let 
\begin{equation}
B_1\subset \mathbb C^{(n+1)(d+1)}
\end{equation}
be the subvariety 
that is equal to 
\begin{equation}
\{0\}\oplus  H^0(\mathcal O_{\mathbf P^1}(d))\oplus \cdots \oplus  H^0(\mathcal O_{\mathbf P^1}(d)).
\end{equation}

Let \begin{equation}\begin{array}{ccc}\tilde \mathbb C^{(n+1)(d+1)} &\stackrel{\pi_1} \rightarrow & 
\mathbb C^{(n+1)(d+1)}\end{array}\end{equation}

 be the blow-up of $\mathbb C^{(n+1)(d+1)}$ along $B_1$. It is clear that
\begin{equation}\begin{array}{cc}
 & \tilde \mathbb C^{(n+1)(d+1)}\\
& \|\\ &
\tilde \mathbb C^{d+1}\times H^0(\mathcal O_{\mathbf P^1}(d))\times\cdots \times H^0(\mathcal O_{\mathbf P^1}(d))\end{array}
\end{equation}
where  $\tilde \mathbb C^{d+1}$ is the blow-up of $\mathbb C^{d+1}$ at the origin.
Now let $\tilde M$ and $\tilde P(\Gamma_{\mathbb L})$ ($\tilde P(\Gamma_{\mathbb L_1})$ ) be the strict transforms of $M$ and 
$P(\Gamma_{\mathbb L})$ ($P(\Gamma_{\mathbb L_1})$ )  respectively. Let
$ E_1$, $E_{1\mathbb L}$ be their exceptional divisors.
 Let $\tilde c_2\in \tilde P(\Gamma_{\mathbb L_1})$ be an inverse of $c_2$ under the map
\begin{equation}\begin{array}{ccc}
 \tilde P(\Gamma_{\mathbb L_1}) &\stackrel{\pi_1|_{ \tilde P(\Gamma_{\mathbb L_1})}} \rightarrow & 
P(\Gamma_{\mathbb L_1}).\end{array} 
\end{equation}
(such $\tilde c_2$ is independent of choice of $q$). 
\bigskip

\begin{lemma} Let $(t_1, \cdots, t_{hd})\in Sym^{hd}(\mathbf P^1)$ be generic. Let
\begin{equation}\begin{array}{cc}  \psi_i=\left|  \begin{array}{cc} f_0(c(t_1)) &  f_2(c(t_1))\\
f_0(c(t_2)) & f_2(c(t_2))\end{array}\right| df_1(c(t_i))+\left|  \begin{array}{cc} f_2(c(t_1)) &  f_1(c(t_1))\\
f_2(c(t_2)) & f_1(c(t_2))\end{array}\right| df_0(c(t_i))& \\  +\left|  \begin{array}{cc} f_1(c(t_1)) &  f_0(c(t_1))\\
f_1(c(t_2)) & f_0(c(t_2))\end{array}\right| df_2(c(t_i))
\end{array}\end{equation}
for $i=3, 4, \cdots, hd+1$. Then
vectors
\begin{equation}\begin{array}{cc} &  \pi_1^\ast (\psi_i), i=3, \cdots, hd+1\\
& \pi_1^\ast (df_0(c(t_1))), \pi_1^\ast (df_0(c(t_2))),\\
& \pi_1^\ast (df_1(c(t_1))), \pi_1^\ast (df_1(c(t_2))),\\
& \pi_1^\ast (df_2(c(t_1))), \pi_1^\ast (df_2(c(t_2)))
\end{array}\end{equation}
are linearly independent in the vector space $(T_{\tilde c}{\tilde{M}})^\ast$, i.e. they form a basis, where 
$\tilde c$ lies in a non-empty open set of  $ \tilde P(\Gamma_{\mathbb L})$, and $\tilde c\neq \tilde c_2$.  \footnote{  The vectors in lemma 2.14 come from the expansion of $\phi_i$: $$\begin{array}{cc}
  \phi_i= \left|  \begin{array}{cc} f_0(c(t_1)) &  f_2(c(t_1))\\
f_0(c(t_2)) & f_2(c(t_2))\end{array}\right| df_1(c(t_i))+\left|  \begin{array}{cc} f_2(c(t_1)) &  f_1(c(t_1))\\
f_2(c(t_2)) & f_1(c(t_2))\end{array}\right| df_0(c(t_i))& \\  +\left|  \begin{array}{cc} f_1(c(t_1)) &  f_0(c(t_1))\\
f_1(c(t_2)) & f_0(c(t_2))\end{array}\right| df_2(c(t_i))
+ \sum_{l=0, j=1}^{l=2, j=2} h_{lj}^i(c) df_l(c(t_j)) & \end{array}$$
where $h_{lj}^i$ are polynomials in $c$.}   

\end{lemma}
\bigskip

\begin{proof} of lemma 2.14: The proof  is long. Thus we divide it into two steps.
It suffices to prove it for special $t_1, \cdots, t_{hd+1}$.  First let
$\tilde c_2$ be decomposed (according to (2.61)) to 
$$\tilde c_2=( \tilde c_2^0(t), \cdots, \tilde c_2^n(t))$$
where $$\tilde c_2^0(t)\in \mathbf P(H^0(\mathcal O_{\mathbf P^1}(d)), \tilde c_2^i(t)\in H^0(\mathcal O_{\mathbf P^1}(d)), i\neq 0.$$
and the zeros of all $\tilde c_2^i, i=0, \cdots, n$ are $\tilde\theta_i^j$.  
\par

{\bf Step 1}: 
 Suppose that
 $$\tilde\theta_i^j, , i=0, \cdots, n, j=1, \cdots, d,$$ are distinct. In this case, we 
let  $t_1, t_2$ be two points on $\mathbf P^1$ satisfying
\begin{equation}
\left|\begin{array}{cc} q|_{\tilde c_2(t_1)}  &  \tilde c_2^{n-1}(t_1) \tilde c_2^n (t_1)\\
q|_{\tilde c_2(t_2)}& \tilde c_2^{n-1}(t_2) \tilde c_2^n (t_2)\end{array}\right|=0.
\end{equation}

Let 
\begin{equation}
f_3=z_0\cdots z_{n-2}( \delta_1 q+\delta_2 z_{n-1} z_n)
\end{equation}
where $$\delta_1=\left|\begin{array}{cc} f_0 (c_2(t_1))  &   f_2(\tilde c_2(t_1))\\
f_0 (c_2(t_2))& f_2(\tilde c_2(t_2))\end{array}\right|, \quad \delta_2=\left|\begin{array}{cc} f_1 (\tilde c_2(t_1))  &   f_0(c_2(t_1))\\
f_1 (\tilde c_2(t_2))& f_0( c_2(t_2))\end{array}\right|.$$
Then let 
$t_3, \cdots, t_{hd}$ be zeros of
\begin{equation}\begin{array}{c} f_3((\tilde c_2(t))\\
\|\\
\tilde c_2^0(t) \cdots \tilde c_2^{n-2}(t) \biggl ( 
\delta_1
q|_{\tilde c_2(t)}+  \delta_2 \tilde c_2^{n-1}(t)\tilde c_2^n(t)\biggr)\\
\| \\
0,\end{array}
\end{equation}
other than $\tilde\theta_0^1, \tilde\theta_1^1$.  Let $t_{hd+1}$ be generic.  
Because $f_0,  q$ are generic,  
the zeros of 
\begin{equation}  \delta_1
q|_{\tilde c_2(t)}+  \delta_2 \tilde c_2^{n-1}(t)\tilde c_2^n(t)
\end{equation}
are distinct and non-zeros. 
Thus \begin{equation}
t_1, t_2, t_3, \cdots, t_{hd+1}
\end{equation}
are distinct.  Let's set-up the coordinates of $\tilde \mathbb C^{(n+1)(d+1)}$.
Let $c_i^j, i=0, \cdots, n, j=0, \cdots, d$ be the coefficients of
$n+1$  tuples of sections of $H^0(\mathcal O_{\mathbf P^1}(d))$. They are the affine coordinates of
$$\mathbb C^{(n+1)(d+1)}.$$
Each section of $\mathcal O_{\mathbf P^1}(d)$  in an analytic neighborhood excluding those with multiple zeros can be
written as 
\begin{equation}
c_i(t)=r_i\Pi_{j=1}^d (t-\theta_i^j).
\end{equation}
for $i\leq n-2, r_i\neq 0$, and
\begin{equation}
c_i(t)=r_i\Pi_{j=1}^d (t-\mu_i^j),
\end{equation}
for $i=n-1, n$ and $r_i\neq 0$. 

Then $r_m, \mu_i^j, \theta_i^j$ are local analytic coordinates for an analytic open set $U_{\tilde \mathbb C^{(n+1)(d+1)}}$ of the blow-up 
\begin{equation}  \tilde \mathbb C^{(n+1)(d+1)}
\end{equation} 
centered around $\tilde c_2\in \tilde \mathbb C^{(n+1)(d+1)}$(now $r_m$ as coordinates of $\tilde \mathbb C^{(n+1)(d+1)}$ could be zeros).  We may assume
$\tilde c_2$ lies in the neighborhood of the coordinates\footnote{If not, we continue to have successive blow-ups till the pre-image of $c_2$ lies 
in the coordinates' neighborhood.} and 
$\tilde c_2$ has specific coordinates
\begin{equation}\begin{array}{cc} &
r_0=0 \\
& \theta_i^j=\tilde \theta_i^j, i=0, \cdots, n-2,  j=1, \cdots, d, \\
& \mu_i^j=\tilde\theta_i^j,  i=n-1, n, j=1, \cdots, d,\\
& r_j=y_j\neq 0, j=1, \cdots, n
\end{array}\end{equation}

Because $t_3, \cdots, t_{hd}$ are distinct and the  last $2d$ of them are non-zeros, by the definition 2.11, in this neighborhood we have another analytic coordinates
\begin{equation} \theta_0^1, \cdots, \theta_{n-2}^d, r_0, \cdots, r_n, \epsilon_1, \cdots, \epsilon_{2d}.
\end{equation}
To see this, we notice that there is an isomorphism $\tilde \mathbb C^{h(d+1)}\to \mathbb C^{h(d+1)}$ outside of $B_1$.  Thus above the coordinates 
(2.76) from definition 2.11 are also local coordinates of $\tilde \mathbb C^{h(d+1)}$.  
In the rest of calculation we use these coordinates for $\tilde \mathbb C^{h(d+1)}$. 
Then the finite set 
\begin{equation} \mathcal B=\{dr_l, d\theta_i^j, d\epsilon_m \}_{l=0, \cdots, n, i\leq n-2, j=1, \cdots, d, m=1, \cdots, 2d}\end{equation}
 is  a basis for $(T_{\tilde c} U_{\tilde \mathbb C^{(n+1)(d+1)}})^\ast$ where
$\tilde c\in U_{\tilde \mathbb C^{(n+1)(d+1)}}$. 
\bigskip

For this step, it is clear that lemma 2.14 follows from 
\bigskip

\begin{lemma}
Let $\mathcal A'$ be the coefficient matrix of $hd+5$ vectors in (2.66) under the basis $\mathcal B$. Then
$\mathcal A'$ has full rank near $\tilde c_2$.

\end{lemma} 

\bigskip

\begin{proof} of lemma 2.15:  
\par
{\it Set-up}: Assume  the blow-up in the lemma 2.15. 
The row vectors in $\mathcal A'$ are vectors in (2.66). We place them from top to bottom in the following order
\begin{equation} \begin{array}{cc} &
\pi_1(\psi_3), \\
& \vdots\\
& \pi_1^\ast (\psi_{hd}),\\
& \pi_1^\ast (\psi_{hd+1})\\
& \pi_1^\ast (df_2(c(t_{1}))),\\
& \pi_1^\ast (df_2(c(t_{2}))),\\
& \pi_1^\ast (df_1(c(t_{1}))),\\
& \pi_1^\ast (df_1(c(t_{2}))),\\
& \pi_1^\ast (df_0(c(t_{1}))),\\
& \pi_1^\ast (df_0(c(t_{2}))).
\end{array}\end{equation}
The basis vectors are listed,  from left to right, as 
$$
d\theta_0^2, \cdots, \widehat d\theta_1^1, \cdots,  d\theta_{n-2}^d, 
d\epsilon_1, \cdots, d\epsilon_{2d}, dr_0, d\theta_0^1, d\theta_1^1, dr_1, \cdots, dr_n.
$$
( $\widehat\cdot $  denotes omitting). So $\mathcal A'$ is a matrix of
size $[hd+5]\times [(n+1)(d+1)]$ where $n\geq 4$. So there are more columns than rows. Let
$\mathcal A$ be the sub square matrix of $\mathcal A'$  which is the coefficient matrix of
vectors in (2.66) under the basis

\begin{equation}\begin{array} {c}
d\theta_0^2, \cdots, \widehat d\theta_1^1, \cdots,  d\theta_{n-2}^d, 
d\epsilon_1, \cdots, d\epsilon_{2d}, \\
dr_0, d\theta_0^1, d\theta_1^1, dr_1, dr_2, dr_{n-1}, dr_n.
\end{array}\end{equation}
This can be arranged because $n\geq 4$. 
Then matrix $\mathcal A$ is the matrix  partitioned as 

\begin{equation}\left(\begin{array}{ccc}
\mathcal A_{1 1} &\mathcal A_{12} &\mathcal A_{13} \\
\mathcal A_{2 1} &\mathcal A_{22} &\mathcal A_{23} \\
\mathcal A_{3 1} &\mathcal A_{32} &\mathcal A_{33}\\
\mathcal A_{4 1} &\mathcal A_{42} &\mathcal A_{43}\\
\mathcal A_{5 1} &\mathcal A_{52} &\mathcal A_{53}
\end{array}\right)\end{equation}
where all $\mathcal A_{ij}$ are matrices of different sizes. We should describe them, one-by-one,  as follows:
\par

In the following $ O(1)$ denotes a polynomial function on $U_{\tilde \mathbb C^{(n+1)(d+1)}}$ of
forms 
\begin{equation}
(t_l-\theta_i^j) O, \epsilon_i O
\end{equation}
where $O$ are polynomial functions on $U_{\tilde \mathbb C^{(n+1)(d+1)}}$. Note that by the definition $$t_l=\tilde\theta_i^j, i\leq n-2$$ are some complex numbers. 
Thus 
$O(1)$ could indicate different functions. \par

(I) $ \mathcal A_{11}$: 
it is a $(hd-2)\times (hd-2)$ matrix. Entries are coefficients of 
$$d\theta_i^j, d\epsilon_m, i\leq n-2, (i, j)\neq (0, 1), (1, 1), m=1, \cdots, 2d$$ for the vectors 
$$  \pi_1^\ast (\psi_i), i=3, \cdots, hd.
$$
For this block matrix we use lemma 2.12 to obtain that 
 its diagonal entries in the order of top-to-bottom are
\begin{equation} 
b_3, \cdots, b_{hd}
\end{equation}
where all $b_i=r_0^2 b_i'$ such that  $b_i'$ are polynomial functions on $U_{\tilde \mathbb C^{(n+1)(d+1)}}$ that do not vanish at $\tilde c_2$.
All the rest of entries are in the form $$r_0^2 O(1),$$ where 
$O(1) $ is defined in (2.81).\par

(II) $\mathcal A_{12}$: it is $(hd-2)\times 1$ matrix. Entries are coefficients of
$dr_0$ for the vectors $$\pi_1^\ast (\psi_i), i=3, \cdots, hd.$$
So from top-to-bottom they are
$$a_{3}, \cdots, a_{hd}$$
where  $a_l=r_0 O(1)$.

(III) $ \mathcal A_{13}$: it is a $(hd-2)\times 6$ matrix. Entries are coefficients of $$d\theta_0^1, d\theta_1^1, dr_1, dr_2, dr_{n-1}, dr_n$$ for the vectors

$$ \pi_1^\ast (\psi_i), i=3, \cdots, hd.$$
So they are all in the form
$$r_0^2 O(1).$$

(IV) $\mathcal A_{21}$: it is a $1\times (hd-2)$ matrix. Entries are coefficients of 
$$d\theta_i^j, d\epsilon_m, i\leq n-2, (i, j)\neq (0, 1), (1, 1), m=1, \cdots, 2d$$   for
the vector 
$$\pi_1^\ast (\psi_{hd+1}).$$
So all entries are in the form $$r_0^2 O$$
 where 
$O $ is a polynomial function on $U_{\tilde \mathbb C^{(n+1)(d+1)}}$. 

\par
(V) $\mathcal A_{22}$: it is a $1\times 1$ matrix. It is the coefficient of $dr_0$ for the vector, 
$$\pi_1^\ast (\psi_{hd+1}).$$
It is in the form
$$r_0 O_{hd+1}$$ where $O_{hd+1}$ is a polynomial function on $U_{\tilde \mathbb C^{(n+1)(d+1)}}$ that does not vanish at $\tilde c_2$.
\par

(VI) $\mathcal A_{23}$: it is a $1\times 6$ matrix. Entries are coefficients of 
$$d\theta_0^1, d\theta_1^1, dr_1, dr_2, dr_{n-1}, dr_n$$ for the vector,
$$\pi_1^\ast (\psi_{hd+1}).$$

So all entries are in the form $$r_0^2 O$$
 where 
$O $ is a polynomial function on $U_{\tilde \mathbb C^{(n+1)(d+1)}}$. 
\par
(VII) $\mathcal A_{31}$: it is a $2\times (hd-2)$ matrix. Entries are coefficients of
$$d\theta_i^j, d\epsilon_m, i\leq n-2, (i, j)\neq (0, 1), (1, 1), m=1, \cdots, 2d$$   for the vectors
$$\pi_1^\ast (df_2(c(t_{1}))),  \pi_1^\ast (df_2(c(t_{2}))).$$
So the entries  are all in the form
$$r_0 O, $$
where 
$O $ is a polynomial function on $U_{\tilde \mathbb C^{(n+1)(d+1)}}$. 

\par

(VIII) $\mathcal A_{32}$: it is $2\times 1$ matrix. Entries are the coefficients of 
$dr_0$ for the vectors
$$\pi_1^\ast (df_2(c(t_{1}))),  \pi_1^\ast (df_2(c(t_{2}))).$$
So they are in
the form of $$ O,$$
where 
$O $ is a polynomial function on $U_{\tilde \mathbb C^{(n+1)(d+1)}}$. 

(IVV) $\mathcal A_{33}$: it is $2\times 6$ matrix. Entries are the coefficients of 
$$d\theta_1^1, d\theta_2^1, dr_1, dr_2, dr_{n-1}, dr_n$$ for the vectors
$$\pi_1^\ast (df_2(c(t_{1}))),  \pi_1^\ast (df_2(c(t_{2}))).$$
So all entries are in the forms of 
$$r_0 O,$$
where 
$O $ is a polynomial function on $U_{\tilde \mathbb C^{(n+1)(d+1)}}$. 

(VV) $\mathcal A_{41}$: it is a $2\times (hd-2)$ matrix. Entries are coefficients of
$$d\theta_i^j, d\epsilon_m, i\leq n-2, (i, j)\neq (0, 1), (1, 1), m=1, \cdots, 2d$$    for
$$
 \pi_1^\ast (df_1(c(t_{1}))),
 \pi_1^\ast (df_1(c(t_{2}))).
$$
All entries are in the form of $$r_0 O,$$
where 
$O $ is a polynomial function on $U_{\tilde \mathbb C^{(n+1)(d+1)}}$. 
\par

(VVI) $\mathcal A_{42}$: it is a $2\times 1$ matrix. Entries are the coefficients of
$dr_0$ for the vector
$$
 \pi_1^\ast (df_1(c(t_{1}))),
 \pi_1^\ast (df_1(c(t_{2}))).
$$
\par
(VVII) $\mathcal A_{43}$: it is a $2\times 6$ matrix. 
Entrices are the coefficients of
$$d\theta_0^1, d\theta_1^1, dr_1, dr_2, dr_{n-1}, dr_n$$ for the vector
$$
 \pi_1^\ast (df_1(c(t_{1}))),
 \pi_1^\ast (df_1(c(t_{2}))).
$$
All entries are in the form 
$$r_0 O,$$
where 
$O $ is a polynomial function on $U_{\tilde \mathbb C^{(n+1)(d+1)}}$.

\par

(VVIII) $\mathcal A_{51}$: it is a $2\times (hd-2)$ matrix. Entries are coefficients of
$$d\theta_i^j, d\epsilon_m, i\leq n-2, (i, j)\neq (0, 1), (1, 1), m=1, \cdots, 2d$$    for
$$
 \pi_1^\ast (df_0(c(t_{1}))),
 \pi_1^\ast (df_0(c(t_{2}))).
$$
\par

(VVIV) $\mathcal A_{52}$: it is a $2\times 1$ matrix. Entries are the coefficients of
$dr_0$ for the vector
$$
 \pi_1^\ast (df_0(c(t_{1}))),
 \pi_1^\ast (df_0(c(t_{2}))).
$$
\par
(VVV) $\mathcal A_{53}$: it is a $2\times 6$ matrix. 
Entries are the coefficients of
$$d\theta_0^1, d\theta_1^1, dr_1, dr_2, dr_{n-1}, dr_n$$ for the vector
$$
 \pi_1^\ast (df_0(c(t_{1}))),
 \pi_1^\ast (df_0(c(t_{2}))).
$$

Next consider the function on $U_{\tilde \mathbb C^{(n+1)(d+1)}}-E_1$ (where $r_0\neq 0$), 
\begin{equation}
\mu(\tilde c)={1\over (r_0)^{2hd+1}}|\mathcal A|
\end{equation}
which is the determinant of
the matrix
\begin{equation} \mathcal A_{r_0}=\left(\begin{array}{ccc}
{1\over r_0^2}\mathcal A_{1 1} &{1\over r_0^2}\mathcal A_{12} &{1\over r_0^2}\mathcal A_{13} \\
{1\over r_0}\mathcal A_{2 1} &{1\over r_0}\mathcal A_{22} &{1\over r_0}\mathcal A_{23}\\
{1\over r_0}\mathcal A_{3 1} &{1\over r_0}\mathcal A_{32} &{1\over r_0}\mathcal A_{33}\\
{1\over r_0}\mathcal A_{4 1} &{1\over r_0}\mathcal A_{42} &{1\over r_0}\mathcal A_{43}\\
\mathcal A_{5 1} &\mathcal A_{52} &\mathcal A_{53}\\
\end{array}\right). \end{equation}
To prove lemma 2.15, it is sufficient to prove the determinant 
$$|\mathcal A_{r_0}|\ at\ r_0=0 $$
is non-zero. 
Unfortunately,  the matrix $\mathcal A_{r_0}$ is not a well-defined at $\tilde c_2$ (i.e. $r_0=0$) 
 because some entries involve  ${(t_l-\theta_i^j)\over r_0}$. But those terms do not 
show up in the computation of its determinant $\mu(\tilde c)$. Hence $\mu(\tilde c)$ can be continuously extended to $\tilde c_2$.
  Let's see this. \smallskip

{\it Computation}: 
As before all expressions in the computation use coordinates in definition 2.11. 
Notice all entries in 

\begin{equation} \left(\begin{array}{c}
{1\over r_0^2}\mathcal A_{1 1} \\
{1\over r_0}\mathcal A_{2 1} \\
{1\over r_0}\mathcal A_{3 1} \\
{1\over r_0}\mathcal A_{4 1}\\
\mathcal A_{5 1} \\
\end{array}\right). \end{equation}
and 
\begin{equation} {1\over r_0^2}\mathcal A_{13} \end{equation}
can be extended to the entire neighborhood $U_{\tilde \mathbb C^{(n+1)(d+1)}}$, and when evaluated at $\tilde c_2$, 
\begin{equation} {1\over r_0^2}\mathcal A_{1 1}
\end{equation}
is a non-zero diagonal matrix (by lemma 2.12).  Also notice that 
\begin{equation} {1\over r_0^2}\mathcal A_{13}|_{\tilde c_2}=0. \end{equation}
Therefore
we can perform the row operations on 
the matrix
$$\mathcal A_{r_0}$$ to reduce
\begin{equation} \left(\begin{array}{c}
{1\over r_0}\mathcal A_{2 1} \\
{1\over r_0}\mathcal A_{3 1} \\
{1\over r_0}\mathcal A_{4 1}\\
\mathcal A_{5 1} \\
\end{array}\right). \end{equation}
to zero matrix.
Hence 
\begin{equation}
\mu(\tilde c)={1\over (r_0)^{2hd+1}}|\mathcal A|
\end{equation}
is a non-zero multiple of 

\begin{equation} \left|\begin{array}{cc}
{1\over r_0}\mathcal A_{22} &{1\over r_0}\mathcal A_{23}\\
{1\over r_0}(\mathcal A_{32}+O_{32}(1)) &{1\over r_0}\mathcal A_{33}\\
{1\over r_0}(\mathcal A_{42}+O_{42}(1)) &{1\over r_0}\mathcal A_{43}\\
\mathcal A_{52}+{1\over r_0}O_{52}(1) &\mathcal A_{53}\\
\end{array}\right|+O(1). \end{equation}
where $O(1), O_{ij}(1)$ represent the determinants  of  matrices,  whose entries are local functions vanishing at $\tilde c_2$. Notice that
in (2.91) the block matrices of first row and first column have sizes $1\times 7$ and $7\times 1$ respectively.  
Therefore by the linearity of determinants,
\begin{equation} \mu(\tilde c)={1\over \rho} \left|\begin{array}{cc}
{1\over r_0}\mathcal A_{22} &{1\over r_0^2}\mathcal A_{23}\\
\mathcal A_{32} &{1\over r_0}\mathcal A_{33}\\
\mathcal A_{42} &{1\over r_0}\mathcal A_{43}\\
r_0\mathcal  A_{52} &\mathcal A_{53}\\
\end{array}\right|+O(1) \end{equation}
for some non-zero complex number $\rho$. 
Now all entries in (2.92) are well-defined functions on $U_{\tilde \mathbb C^{(n+1)(d+1)}}$. 
Therefore it suffices to show the non-degeneracy of this $7\times 7$ matrix 
\begin{equation} \left.\left(\begin{array}{cc}
{1\over r_0}\mathcal A_{22} &{1\over r_0^2}\mathcal A_{23}\\
\mathcal A_{32} &{1\over r_0}\mathcal A_{33}\\
\mathcal A_{42} &{1\over r_0}\mathcal A_{43}\\
r_0\mathcal  A_{52} &\mathcal A_{53}\\
\end{array}\right)\right|_{\tilde c_2}
\end{equation}
at the point $\tilde c_2$. Because $t_{hd+1}$ is generic, also $q$ is generic and $\tilde\theta_i^j$ are distinct, , 
it suffices to prove $6\times 6$ matrix
\begin{equation} \left.\left(\begin{array}{c}
{1\over r_0}\mathcal A_{33}\\
{1\over r_0}\mathcal A_{43}\\
\mathcal A_{53}
\end{array}\right)\right|_{\tilde c_2}
\end{equation}
 is non-degenerate.  

Let $\lambda_1$ be the determinant of (2.94).  Also let
$$ g_i(t)=({1\over r_0}\pi_1^\ast ({ z_i\partial  f_1\over \partial z_i}))|_{\tilde c_2(t)}.$$

Using the coordinates in definition 2.11, we compute that this determinant
\begin{equation} \lambda_1=\left. det \left (\begin{array}{c}
{1\over r_0}\mathcal A_{33}\\
{1\over r_0}\mathcal A_{43}\\
\mathcal A_{53}
\end{array}\right) \right|_{\tilde c_2}
\end{equation}
is equal to 

\begin{equation} \lambda_2 \left|\begin{array}{cccccc}
{1\over t_1-\tilde\theta_0^1} &{1\over t_1-\tilde\theta_1^1} &1 &1  & 1 & 1  \\
{1\over t_2-\tilde\theta_0^1} &{1\over t_2-\tilde\theta_1^1} &1 &1 & 1 & 1 \\
{\partial \pi_1^\ast( f_1)(\tilde c_2 (t_1))\over r_0\partial \theta_0^1} &{\partial \pi_1^\ast( f_1)(\tilde c_2 (t_1))\over r_0\partial \theta_1^1} &
g_1(t_1) & g_2(t_1) &
 g_{n-1}(t_1)& g_n(t_1) \\
{\partial \pi_1^\ast( f_1)(\tilde c_2 (t_2))\over r_0\partial \theta_0^1} &{\partial \pi_1^\ast( f_1)(\tilde c_2 (t_2))\over r_0\partial \theta_1^1} &
g_1(t_2) & g_2(t_2) &
 g_{n-1}(t_2)& g_n(t_2) \\
 {\partial \pi_1^\ast( f_0)(\tilde c_2 (t_1))\over \partial \theta_0^1} & {\partial \pi_1^\ast( f_0)(\tilde c_2 (t_1))\over \partial \theta_1^1}
& (z_1{\partial f_0 \over \partial z_1})|_{\tilde c_2(t_1)} &(z_2{\partial f_0 \over \partial z_2})|_{\tilde c_2(t_1)} &
 (z_{n-1}{\partial f_0 \over \partial z_{n-1}})|_{\tilde c_2(t_1)}& (z_n{\partial f_0 \over \partial z_n})|_{\tilde c_2(t_1)}\\
{\partial \pi_1^\ast( f_0)(\tilde c_2 (t_2))\over \partial \theta_0^1} & {\partial \pi_1^\ast( f_0)(\tilde c_2 (t_2))\over \partial \theta_1^1}
& (z_1{\partial f_0 \over \partial z_1})|_{\tilde c_2(t_2)} &(z_2{\partial f_0 \over \partial z_2})|_{\tilde c_2(t_2)} & 
(z_{n-1}{\partial f_0 \over \partial z_{n-1}})|_{\tilde c_2(t_2)}& (z_n{\partial f_0 \over \partial z_n})|_{\tilde c_2(t_2)}
\end{array}\right|. \end{equation}

where 
\begin{equation}
\lambda_2=\prod_{i=1}^2 ({\pi_1^\ast(f_2)\over y_1y_2y_{n-1}y_nr_0})|_{\tilde c_2(t_i)}\neq 0. 
\end{equation}

Because $q$ is generic, the two  middle row vectors of (2.96)  

\begin{equation}\begin{array}{cccccc}
({\partial \pi_1^\ast( f_1)(\tilde c_2(t_1))\over r_0\partial \theta_0^1} &{\partial \pi_1^\ast( f_1)(\tilde c_2(t_1))\over r_0\partial \theta_1^1} &
g_1(t_1) & g_2(t_1) &
 g_{n-1}(t_1)& g_n(t_1) )\\
({\partial \pi_1^\ast( f_1)( \tilde c_2(t_2))\over r_0\partial \theta_0^1} &{\partial \pi_1^\ast( f_1)(\tilde c_2(t_2))\over r_0\partial \theta_1^1} &
g_1(t_2) & g_2(t_2) &
 g_{n-1}(t_2)& g_n(t_2) )
\end{array}\end{equation}
evaluated at $\tilde c_2$ are generic vectors in $\mathbb C^6$. Thus 
to show (2.96) is non-zero, it suffices to show that the minor

\begin{equation} \begin{array}{c} Jac\\
\|\\
\left|\begin{array}{cccccc}
{1\over t_1-\tilde\theta_0^1}  &1  & 1 & 1  \\
{1\over t_2-\tilde\theta_0^1} &1 & 1 & 1 \\
 {\partial \pi_1^\ast( f_0)(\tilde c_2 (t_1))\over \partial \theta_0^1}  & (z_2{\partial f_0 \over \partial z_2})|_{\tilde c_2(t_1)} &
 (z_{n-1}{\partial f_0 \over \partial z_{n-1}})|_{\tilde c_2(t_1)}& (z_n{\partial f_0 \over \partial z_n})|_{\tilde c_2(t_1)}\\
{\partial \pi_1^\ast( f_0)(\tilde c_2 (t_2))\over \partial \theta_0^1} & (z_2{\partial f_0 \over \partial z_2})|_{\tilde c_2(t_2)} &
 (z_{n-1}{\partial f_0 \over \partial z_{n-1}})|_{\tilde c_2(t_2)}& (z_n{\partial f_0 \over \partial z_n})|_{\tilde c_2(t_2)}
\end{array}\right|\\
\nparallel\\
0. \end{array}\end{equation}
(Remove the 2nd and 3rd columns in (2.96)).
We further compute to have
\begin{equation}\begin{array}{c}  Jac\\
\|\\
({1\over t_1-\tilde\theta_0^1}-{1\over t_2-\tilde\theta_0^1} )\left|\begin{array}{ccc}
1&1&1\\  (z_2{\partial  f_0)\over \partial z_2})|_{c_2(t_1)} &
(z_{n-1}{\partial  f_0)\over \partial z_{n-1}})|_{c_2(t_1)} &(z_n{\partial  f_0)\over \partial z_n})|_{c_2(t_1)} \\
(z_2{\partial  f_0)\over \partial z_2})|_{c_2(t_2)} &
(z_{n-1}{\partial  f_0)\over \partial z_{n-1}})|_{c_2(t_2)} &(z_n{\partial  f_0)\over \partial z_n})|_{c_2(t_2)} 
\end{array}\right|. \end{array}
\end{equation}

Suppose $Jac=0$. Because $q$ is generic, $(t_1, t_2)$ must a generic point in $\mathbb C^2$.  Hence  
$$c_2(t)=\{[0, c_2^1(t), \cdots, c_2^{n-1}(t), c_2^n(t)]\}$$
lies in the hypersurface  $$ z_2{\partial  f_0\over \partial z_2}+z_{n-1}{\partial  f_0\over \partial z_{n-1}}+z_n{\partial  f_0\over \partial z_n}=0.$$
Similarly we can conclude that $\{c_2(t)\}=C_2$ (image of $\mathbf P^1$) lies in the hypersurface
$$H_j=\{\Sigma_{i\neq j } z_i{\partial f_0\over \partial z_i}=0\}$$
for $j=1, \cdots, n$.  
Since $f_0$ is generic,$$\cap_j H_j$$ is finite. Then $C_2$ has to be a point.  This contradicts our assumption for the step 1.  

 Hence 
\begin{equation} \left|\begin{array}{ccc}
1&1&1\\  (z_2{\partial  f_0)\over \partial z_2})|_{c_2(t_1)} &
(z_{n-1}{\partial  f_0)\over \partial z_{n-1}})|_{c_2(t_1)} &(z_n{\partial  f_0)\over \partial z_n})|_{c_2(t_1)} \\
(z_2{\partial  f_0)\over \partial z_2})|_{c_2(t_2)} &
(z_{n-1}{\partial  f_0)\over \partial z_{n-1}})|_{c_2(t_2)} &(z_4{\partial  f_0)\over \partial z_4})|_{c_2(t_2)} 
\end{array}\right|\neq 0. 
\end{equation}
Therefore 
$$Jac\neq 0.$$
This completes the proof of lemma 2.15, thus the first step of lemma 2.14.

\end{proof}

{\bf Step 2}:   The following formulation of  step 2  can be viewed as  the generalization of step 1. Both steps have the same formulation after the blow-ups. \par
Let's assume $\tilde \theta_i^j$ are not distinct.  
 Notice $\omega$ depends on the choice of $t_i, i=1, \cdots, hd+1$.
The main  idea we would like to get across in this paper is that for generic choice of $t_i$, $\omega$ is non-zero near $c_2$ (excluding $c_2$) .
 However it seems to be completely helpless in the computation of $\omega$ for generic $t_i$. This situation changes if we use a special choice of $t_i$ such as in step 1, in which case 
the Jacobian matrix in $\omega$ breaks down to manageable block matrices in (2.80).  The choice of such $t_i$ requires that
$\tilde \theta_i^j$ are distinct. If they aren't (i.e. $\tilde c_2$ has  multiple intersection points with coordinates' plane), 
 the step 1 will fail because the main block in (2.80) becomes degenerate. 
So to resolve this,  we will use successive blow-ups in the following to reduce the multiple $\tilde \theta_i^j$ to distinct $\tilde \theta_i^j$. 
The pull-back of $\omega$ by the blow-ups for generic $t_i$ is exactly the same as (5.10) which is that in step 1. 

\bigskip

Let's discuss the blow-ups.  They are all analytic which means the total spaces and the centers are analytic in an analytic neighborhood. 
Suppose we already have the blow-up $\tilde \mathbb C^{h(d+1)}$ as in (2.62). 
 Recall  $\tilde c_2^i\in H^0(\mathcal O_{\mathbf P^1}(d))$ are  sections such that 
$$ div(\tilde c_2^i)=\sum_{j=1}^{d}\tilde \theta_i^j. $$
for $i=0, \cdots, n$.
We may assume
multiple zeros are
$$\tilde\theta_\alpha^{\beta}$$
among all $\tilde\theta_i^j$, 
where $\alpha, \beta$ are finite numbers less than $h$ and $d+1$ respectively. Let 
$$m(\tilde c_2)$$ be the number of the pairs $(\alpha, \beta)$ which is the number
of multiple zeros $\theta_i^j$.   The worst case is 
$m(\tilde c_2)=hd$, in which case $\tilde c_2$ represents a constant map $\mathbf P^1\to \mathbf P^n$. This multiplicity is well-defined for the map $\tilde c_2$, but depends on the coordinates $z_i$. 
Next we define the successive blow-ups which are controlled by the multiplicity $m(\tilde c_2)$. \par
{\it First blow-up}.  First we should note the coordinates $r_i, \theta_i^j$ are well-defined analytic coordinates for $\tilde M$ around the point
$\tilde c_2$, even if $\theta_i^j$ are not distinct.\footnote{ This is actually true in a Zariski topology. The analytic blow-ups we are using here can be replaced by 
the algebraic ones.} We define  $B_2$ to be the  analytic subvariety 
$$
\{ \tilde\theta_\alpha^{\beta}-\theta_\alpha^\beta=0, for \ all\ \alpha, \beta\}. $$
 of an analytic open set $U_{\tilde \mathbb C^{h(d+1)}}$ of $\tilde \mathbb C^{h(d+1)}$, centered around $\tilde c_2$. 
For the simplicity we denote $U_{\tilde \mathbb C^{h(d+1)}}$ by $\tilde \mathbb C^{h(d+1)}$.
We blow-up $\tilde \mathbb C^{h(d+1)}$ along $B_2$ to obtain the first blow-up map
$$\begin{array}{ccc}
N & \stackrel{\pi_2} \rightarrow & \tilde \mathbb C^{h(d+1)}.
\end{array}$$
Let $\tilde c_2^{(1)}$ be the inverse  of $\tilde c_2$ in the strict transform of $\tilde P(\Gamma_{\mathbb L_1})$.
Using the coordinates $\theta_i^j, r_i$ above,  we can find an open set $U_N$ centered around $\tilde c_2^{(1)}$ such that
$U_N$ is analytically isomorphic to 
\begin{equation} U_N\simeq  \mathbb C^{d+1}\times \cdots \times \mathbb C^{d+1},\end{equation}
 where each copy $\mathbb C^{d+1}$ corresponds to $\theta_i^j$ with the same $i$. 
Precisely the i-th copy  $\mathbb C^{d+1}$ in (2.102) is determined by the defining equations
$$\pi_2^\ast(\theta_i^j)=0, \pi_2^\ast (r_i)=0.$$
Next fix the isomorphism, 
\begin{equation}\begin{array}{ccc} \mathbb C^{d+1} &  \simeq  & H^0(\mathcal O_{\mathbf P^1}(d))\\
(r_i, (\theta_i^0)^{(1)}, \cdots, (\theta_i^d)^{(1)})& \rightarrow &  r_i \sum_{j=0}^d(t-(\theta_i^j)^{(1)}). \end{array}\end{equation}
Then we obtain the new evaluation map $Z^{(1)}$, 
\begin{equation}\begin{array}{ccc}
U_N\times \mathbf P^1 &\stackrel{Z^{(1)}} \rightarrow & \mathbf P^n \\
(\tilde c^{(1)}, t) &\rightarrow & [\tilde c^{(1)}(t)]
.\end{array}\end{equation}
The map $Z^{(1)}$  turns $U_N$ into a family of maps
\footnote{ Rational curves in this family do not lie on hypersurfaces in the family $\mathbb L$. So they are  not from the family 
$\tilde P(\Gamma_{\mathbb L})$.  But their  transformations outside of the exceptional divisor of the blow-up do. This transformation is nothing else but the isomorphism of the blow-up outside of the exceptional divisor through the identification (2.103). Then $\tilde c_2^{(1)}$ is a specialization of the generic member in this family $U_N$ of rational curves. }, $$\mathbf P^1\to \mathbf P^n.$$
Using (2.102) and (2.103), we let
$$
(\tilde c_2)^{(1)}= ((\tilde c_2)_0^{(1)}(t), \cdots, (\tilde c_2)_n^{(1)}(t)),
$$

{\it Continuing blow-ups}.  If  $\tilde c_2^{(1)}(t)$ still has multiple zeros with the planes $\{z_i=0\}$, 
i.e. $(\tilde c_2)_i^{(1)}(t)=0$ for all $i$ have multiple zeros or equivalently 
$$m(\tilde c_2^{(1)})>1, $$ we continue the blow-ups along the multiple zeros as above. 
This type of  blow-ups can continue. 
The last blow-up is obtained when the multiplicity $m(\tilde c_2^{(\kappa)})$ is reduced to $1$ ( need to repeat (2.102) and (2.103) in each blow-up to have
the well-defined $m(\tilde c_2^{(\kappa)})$ ). 
To see that 
 $m(\tilde c_2)$ will be reduced to 1 (i.e. the blow-ups will stop), we go back to the first blow-up. 
Let $(P(\Gamma_{\mathbb L_1}))^{(1)}$ be the strict transform of $\tilde P(\Gamma_{\mathbb L_1})$ under the first blow-up, and
$\tilde c_2^{(1)}$ be a chosen inverse of $\tilde c_2$ in $(P(\Gamma_{\mathbb L_1}))^{(1)}$.  Since $\tilde P(\Gamma_{\mathbb L_1})$
does not lie in the  $B_2$, the exceptional divisor $\mathcal D$ of $(P(\Gamma_{\mathbb L_1}))^{(1)}$ does not lie in 
	the tangent bundle of $B_2$ ( $\mathcal D$ is a subvariety in the projectivization of the normal bundle of $B_2$). 
Therefore  $\tilde c_2^{(1)}$ when regarded as a map  $\mathbf P^1\to \mathbf P^4$  has multiplicity
$$m(\tilde c_2^{(1)})$$ strictly  less than $m(\tilde c_2)$.
After such successive blow-ups, we obtain the inverse $\tilde c_2^{(\kappa)}$ of $\tilde c_2$ whose multiplicity
$$m(\tilde c_2^{(\kappa)})$$ is $1$. 
We let $\pi_3$ be the composition of all such blow-ups. 
So we obtain the birational map

$$ \begin{array}{ccc}
(\tilde \mathbb C^{h(d+1)})^{(\kappa)} &\stackrel{\pi_3}\rightarrow & \tilde \mathbb C^{h(d+1)}.
\end{array}$$
As above, there are a determined open set $U^{(\kappa)}$ of  $(\tilde \mathbb C^{h(d+1)})^{(\kappa)}$, centered at $\tilde c_2^{(\kappa)}$, and
a determined  analytic isomorphism 
\begin{equation}
U^{(\kappa)}\simeq (\mathbb C^{d+1})^{\oplus h}.
\end{equation}

We denote the natural affine coordinates of $$(\mathbb C^{d+1})^{\oplus 5}$$
in (2.105)  by $c^{(\kappa)}$.

\bigskip

{\it Computation of the pull-back by blow-ups}. 
Next we apply the successive blow-ups $\pi_3$ to prove the lemma. 
We would like to show that after the blow-ups we have the same set-up as for the step 1. 
Let $\tilde P(\Gamma_{\mathbb L})^{(\kappa)}$ be the strict transform of $\tilde P(\Gamma_{\mathbb L})$ under $\pi_3$. 
  Then $U^{(\kappa)}\cap \tilde P(\Gamma_{\mathbb L})^{(\kappa)}$ parametrizes a family of rational curves in $\mathbf P^{4}$
(these rational curves do not lie on the hypersurfaces in $\mathbb L$).

These successive blow-ups can be stated in the following diagram

\begin{equation}\begin{array}{ccccccc}
\tilde P(\Gamma_{\mathbb L})^{(\kappa)}& \stackrel{P} \leftarrow & W &\subset  &(\tilde \mathbb C^{h(d+1)})^{(\kappa)}\times \mathbb L &\stackrel{Proj.}  \rightarrow & \mathbb L\\
\downarrow\scriptstyle{ \pi_3} & & \downarrow&& \downarrow\scriptstyle{ \pi_3\times {identity}}  && ||\\
\tilde P(\Gamma_{\mathbb L})  &\stackrel{P} \leftarrow & \Gamma_{\mathbb L} &\subset  & \tilde \mathbb C^{h(d+1)}\times \mathbb L  & \stackrel{Proj.}  \rightarrow & \mathbb L
\end{array}\end{equation}
where $W$ is the strict transform of $\Gamma_{\mathbb L}$ and $ \pi_3$ is birational.  

Restrict this diagram to a small neighborhood $U^{(\kappa)} $ of $ \tilde c_2^{(\kappa)}$, we have the map $\pi_3$
\begin{equation}\begin{array}{ccc}
U^{(\kappa)} &\rightarrow &\tilde \mathbb C^{h(d+1)}\\
\tilde c^{(\kappa) }&\rightarrow &  \tilde c.\end{array}\end{equation}
Notice the map $\tilde c^{(\kappa)} \to   \tilde c$ is the restriction of $\pi_3$ which is an isomorphism outside
of $(\pi_1\circ \pi_3)^{-1}(B_2)$.  This diagram does not explain when $\pi_3$ stops. The termination of the successive blow-ups depends on the 
higher orders of $\tilde P(\Gamma_{\mathbb L})^{(\kappa)}$ explained above. In addition
 blow-ups themselves also rely on the natural coordinates of $M$. 
Next we explain the choice of our new $t_i$ which are different from those in step 1 because of the multiplicity.  Let
$$f_3=z_0\cdots z_{n-2}( \delta_1 q+\delta_2 z_{n-1} z_n),
$$
and $t_1, t_2$ be generic numbers satisfying
$$
\left|\begin{array}{cc} q|_{\tilde c_2^{(\kappa)}(t_1)}  &  (\tilde c_2^{n-1})^{(\kappa)}(t_1) (\tilde c_2^n)^{(\kappa)}(t_1)\\
q|_{\tilde c_2^{(\kappa)} (t_2)}& (\tilde c_2^{n-1})^{(\kappa)}(t_2) (\tilde c_2^n)^{(\kappa)} (t_2))\end{array}\right|=0.
$$
where $\tilde c_2^{(\kappa)}=( (\tilde c_2^0)^{(\kappa)}, \cdots,  (\tilde c_2^n)^{(\kappa)})$, and
$$\delta_1=\left|\begin{array}{cc} f_0(\tilde c_2^{(\kappa)}(t_1))  &   f_2(\tilde c_2^{(\kappa)}(t_1))\\
f_0 (\tilde c_2^{(\kappa)} (t_2))& f_2(\tilde c_2^{(\kappa)}(t_2))\end{array}\right|, \quad \delta_2=\left|\begin{array}{cc} f_1 (\tilde c_2^{(\kappa)}(t_1))  &   f_0(\tilde c_2^{(\kappa)}(t_1))\\
f_1 (\tilde c_2^{(\kappa)}(t_2))& f_0( \tilde c_2^{(\kappa)})(t_2))\end{array}\right|.$$
Let $(\tilde \theta_i^j)^{(\kappa)}$ be the zeros of 
$(\tilde c_2^i)^{(\kappa)}(t)=0, i=0, \cdots, n$. 
Then let 
$t_3, \cdots, t_{hd}$ be zeros of
\begin{equation}\begin{array}{c} f_3(z)|_{\tilde c_2^{(\kappa)}(t)}\\
\|\\
(\tilde c_2^0)^{(\kappa)} (t) \cdots(\tilde c_2^{n-2})^{(\kappa)}(t) \biggl ( 
\delta_1
q|_{\tilde c_2^{(\kappa)}(t)}+  \delta_2 (\tilde c_2^{n-1})^{(\kappa)}(t) (\tilde c_2^n)^{(\kappa)}(t)\biggr)\\
\|\\
0.\end{array}
\end{equation}
other than $(\tilde \theta_0^1)^{(\kappa)}, (\tilde \theta_1^1)^{(\kappa)}$.  Let $t_{hd+1}$ be generic.  

With such a special choice of $t_i$, 
we can repeat the step 1 for  the same quintics $f_i(z)$ to construct
$ \omega^{(\kappa)}$ where $\omega^{(\kappa)}$ has the current choice of $t_i$ and the variables are $\tilde c^{(\kappa)}$ for $U^{(\kappa)}$. Specifically, we define
\begin{equation} \phi_i^{(\kappa)}= d
\left|  \begin{array}{ccc} f_2(\tilde c^{(\kappa)}(t_i)) & f_1(\tilde c^{(\kappa)}(t_i)) & f_0(\tilde c^{(\kappa)}(t_i))\\
f_2(\tilde c^{(\kappa)}(t_1)) & f_1(\tilde c^{(\kappa)}(t_1)) & f_0(\tilde c^{(\kappa)}(t_1))\\
f_2(\tilde c^{(\kappa)}(t_2)) & f_1(\tilde c^{(\kappa)}(t_2)) & f_0(\tilde c^{(\kappa)}(t_2))
\end{array}\right|\end{equation}
for $i=3, \cdots, hd+1$
\begin{equation}
 \omega^{(\kappa)}=\wedge_{i=3}^{hd+1} \phi_i^{(\kappa)}\in H^0(\Omega_{U^{(\kappa)}}).
\end{equation}
Then we repeat the same process (but without blow-up $\pi_1$) in step 1 to construct the Jacobian 
matrix 
$$\mathcal A^{(\kappa)}$$ 
with respect to $\omega^{(\kappa)}$ and the point $\tilde c_2^{(\kappa)}$.
Then we use the same process in step 1 to calculate $det( \mathcal A^{(\kappa)})$. We obtain 
same result as in lemma 5.4,  i.e.
\begin{equation}
det( \mathcal A^{(\kappa)})\neq 0. 
\end{equation}
at the origin $\tilde c_2^{(\kappa)}$.  Next we discuss the relation between 
$det(\mathcal A)$ and $det( \mathcal A^{(\kappa)}) $.
We would like see that,  as in step 1
$det(\mathcal A)$ is for the Jacobian of the generic points of $\tilde P(\Gamma_L)$, 
 $det( \mathcal A^{(\kappa)}) $ now is for the Jacobian of the generic points of $(\tilde P(\Gamma_L))^{(\kappa)}$.
Notice that $\pi_3$ is an isomorphism outside of exceptional divisor $\pi_3^{-1}(B_2)$. 
So the equations
$$\left|  \begin{array}{ccc} f_2(\tilde c^{(\kappa)}(t_i)) & f_1(\tilde c^{(\kappa)}(t_i)) & f_0(\tilde c^{(\kappa)}(t_i))\\
f_2(\tilde c^{(\kappa)}(t_1)) & f_1(\tilde c^{(\kappa)}(t_1)) & f_0(\tilde c^{(\kappa)}(t_1))\\
f_2(\tilde c^{(\kappa)}(t_2)) & f_1(\tilde c^{(\kappa)}(t_2)) & f_0(\tilde c^{(\kappa)}(t_2))
\end{array}\right|=0$$
$i=3, \cdots, hd+1$ defines the scheme 
$$\tilde P(\Gamma_{\mathbb L})^{(\kappa)}=\pi_3^{-1} (\tilde P(\Gamma_{\mathbb L}))$$
in $U^{(\kappa)}-\pi_3^{-1}(B_2)$. This shows that 
  on the open set 
\begin{equation}  U^{(\kappa)}\cap \tilde P(\Gamma_{\mathbb L})^{(\kappa)}-\pi_3^{-1}(B_2) \end{equation}

\begin{equation} (\pi_3)^\ast (det(\mathcal A))=g \cdot det( \mathcal A^{(\kappa)}) \end{equation}
where $g$ is a function nowhere zero on  
$$ \tilde P(\Gamma_{\mathbb L})^{(\kappa)}-\pi_3^{-1}(B_2).$$
Because of (2.111), $det( \mathcal A)\neq 0$ at a point in 
$$ \pi_1\circ \pi_3(U^{(\kappa)})\cap P(\Gamma_{\mathbb L})-\{c_2\}. $$

This proves lemma 2.14.

\end{proof}

\bigskip

\begin{proof}  of lemma 2.13: 
We first calculate 1-form $\phi_i$ on $M$ evaluated at general $c_g\in M$. 
For $i=3, \cdots, hd+1$, 
\begin{equation}\begin{array}{cc}
  \phi_i= d
\left|  \begin{array}{ccc} f_2(c(t_i)) & f_1(c(t_i)) & f_0(c(t_i))\\
f_2(c(t_1)) & f_1(c(t_1)) & f_0(c(t_1))\\
f_2(c(t_2)) & f_1(c(t_2)) & f_0(c(t_2))
 \end{array}\right|&\\ 
=\left|  \begin{array}{cc} f_0(c_g(t_1)) &  f_2(c_g(t_1))\\
f_0(c_g(t_2)) & f_2(c_g(t_2))\end{array}\right| df_1(c(t_i))+\left|  \begin{array}{cc} f_2(c_g(t_1)) &  f_1(c_g(t_1))\\
f_2(c_g(t_2)) & f_1(c_g(t_2))\end{array}\right| df_0(c(t_i))& \\  +\left|  \begin{array}{cc} f_1(c_g(t_1)) &  f_0(c_g(t_1))\\
f_1(c_g(t_2)) & f_0(c_g(t_2))\end{array}\right| df_2(c(t_i))
+ \sum_{l=0, j=1}^{l=2, j=2} h_{lj}^i(c_g) df_l(c(t_j)) & \end{array}\end{equation}

where $h_{lj}^i$ are polynomials in $c$. \par

The lemma 2.14 implies that 
for a generic choice of $$f_0, f_1, f_2, c_g, t_1, \cdots, t_{hd+1},$$

\begin{equation}\begin{array}{cc} & 
\left|  \begin{array}{cc} f_0(c_g(t_1)) &  f_2(c_g(t_1))\\
f_0(c_g(t_2)) & f_2(c_g(t_2))\end{array}\right| df_1(c(t_i))+\left|  \begin{array}{cc} f_2(c_g(t_1)) &  f_1(c_g(t_1))\\
f_2(c_g(t_2)) & f_1(c_g(t_2))\end{array}\right| df_0(c(t_i))\\ &
+\left|  \begin{array}{cc} f_1(c_g(t_1)) &  f_0(c_g(t_1))\\
f_1(c_g(t_2)) & f_0(c_g(t_2))\end{array}\right| df_2(c(t_i)) \end{array}\end{equation} for
$i=3, \cdots, hd+1$, and
\begin{equation} df_l(c(t_j)), l=0, 1, 2, j=1, 2
\end{equation}
are $hd+5$ linearly independent vectors in $(T_{c_g}M)^\ast$ for a generic $c_g\in P(\Gamma_{\mathbb L})$ (not at special point $c_2$ ).
 i.e. they
form a basis of the vector space $(T_{c_g}M)^\ast$.

This implies the set of 1-forms $\{\phi_i\}_{ i=3, \cdots, hd+1}$ is a linearly independent set in $(T_{c_g}M)^\ast$ for
generic $c_g\in (T_{c_g}M)^\ast$.  Thus $\omega$ is nowhere zero when
it is evaluated at  generic points of $P(\Gamma_{\mathbb L})$.  The lemma 2.13 is proved.

\end{proof}

\bigskip

\subsubsection{Ranks of differential sheaves}
\quad\smallskip

\bigskip

\begin{proof} of proposition 1.3: 
Let $\mathcal N$ be the submodule of global sections, $H^0(\Omega_{M})$
generated by
elements
\begin{equation} \phi_i= d
\left|  \begin{array}{ccc} f_2(c(t_1)) & f_1(c(t_1)) & f_0(c(t_1))\\
f_2(c(t_2)) & f_1(c(t_2)) & f_0(c(t_2))\\
f_2(c(t_i)) & f_1(c(t_i)) & f_0(c(t_i))\end{array}\right|\end{equation}
for $i=3, \cdots, hd+1$. 
Recall that 
$$\left|  \begin{array}{ccc} f_2(c(t_1)) & f_1(c(t_1)) & f_0(c(t_1))\\
f_2(c(t_2)) & f_1(c(t_2)) & f_0(c(t_2))\\
f_2(c(t_i)) & f_1(c(t_i)) & f_0(c(t_i))\end{array}\right|=0,$$
for $i=3, \cdots, hd+1$
define the scheme $P(\Gamma_{\mathbb L})$ for a small $\mathbb L$.  Then
 \begin{equation}
\widetilde {({H^0(\Omega_{M})\over \mathcal N})}\otimes \mathcal O_{P(\Gamma_{\mathbb L})}
\simeq \Omega_{P(\Gamma_{\mathbb L})},
\end{equation}
where $\widetilde {(\cdot)}$ denotes the sheaf associated to the module $(\cdot)$.

Therefore

\begin{equation}\begin{array}{cc} &
 {({H^0(\Omega_{M})\otimes k(c_g)\over \mathcal N\otimes k(c_g)})}
\simeq \Omega_{P(\Gamma_{\mathbb L})}\otimes k(c_g)\\&
=(\Omega_{P(\Gamma_{\mathbb L})})|_{(\{c_g\})},\end{array}\end{equation}
where $k(c_g)=\mathbb C$ is the residue field at generic $$c_g\in P(\Gamma_{\mathbb L}).$$ 
Notice two sides of (2.119) are finitely dimensional  linear spaces over $\mathbb C$.
\begin{equation}\begin{array}{cc} &
dim_{\mathbb C} ((\Omega_{P(\Gamma_{\mathbb L})})|_{(\{c_g\})})\\
&=dim_{\mathbb C}(H^0(\Omega_{M})\otimes k(c_g))
- dim(\mathcal N\otimes k(c_g))\end{array}
\end{equation}

Since \begin{equation}
dim_{\mathbb C} ((\Omega_{P(\Gamma_{\mathbb L})})|_{(\{c_g\})}))=dim(T_{c_g}P(\Gamma_{\mathbb L}))
\end{equation}

\begin{equation}
dim(T_{c_g}P(\Gamma_{\mathbb L}))=dim(M) -dim(\mathcal N\otimes k(c_g)).
\end{equation} 

By lemma 2.13, 
$$dim(\mathcal N\otimes k(c_g))=deg(\omega)=hd-1.$$

The proposition 1.3 is proved. 

\end{proof}

\bigskip

\begin{proof} of theorem 1.1: 
 So 
by  proposition 1.3,  
\begin{equation}
(n+1)(d+1)-(hd-1)
\end{equation}
is the dimension of the Zariski tangent space $T_{c_0}P(\Gamma_{\mathbb L})$. 
Furthermore using lemma 2.8  par (b) and  lemma 2.9, we obtain that
\begin{equation}
T_{c_0}\Gamma_{f_0}=(n+1)(d+1)-(hd-1)-2.
\end{equation} 
Again  using lemma 2.8 part (a) ,  we obtain that
\begin{equation}
dim(H^0(N_{c_0/X_0})) =(n+1)(d+1)-hd-5.
\end{equation}

Now we consider the exact sequence of sheaf modules on $\mathbf P^1$, 
\begin{equation}\begin{array}{ccccccccc}
0 &\rightarrow & N_{c_0/X_0} &\rightarrow &  N_{c_0/\mathbf P^n}&\rightarrow &
c_0^\ast(N_{X_0/\mathbf P^n}) &\rightarrow & 0.\end{array}\end{equation}
This induces the exact sequence of finite dimensional linear spaces
\begin{equation}\begin{array}{ccccccccccc}
0 &\rightarrow & H^0(N_{c_0/X_0}) &\rightarrow & H^0( N_{c_0/\mathbf P^n}) &\rightarrow &
H^0(c_0^\ast(N_{X_0/\mathbf P^n})) &\rightarrow & H^1(N_{c_0/X_0}) &\rightarrow & 0.\end{array}\end{equation}

Using  Eurler sequence for $\mathbf P^n$, we obtain
\begin{equation} dim(H^0( N_{c_0/\mathbf P^n}))=(n+1)(d+1)-4.\end{equation} 
Using adjunction formula, we obtain  \begin{equation}  
 dim(H^0(c_0^\ast(N_{X_0/\mathbf P^n}))=hd+1.\end{equation}

Then by the calculations in (2.125), (2.128) and (2.129), we obtain 

\begin{equation} 
dim(H^0 (c_0^\ast(N_{X_0/\mathbf P^n}))=dim(H^0( N_{c_0/\mathbf P^n}))-dim(H^0(N_{c_0/X_0})) ).
\end{equation}

Hence by (2.127)
\begin{equation} dim(H^1(N_{c_0/X_0}))=0.
\end{equation}

Theorem 1.1 in the Calabi-Yau case is proved.

\end{proof}

\bigskip

\subsection{Hypersurface of general type}\smallskip\quad

In this subsection, we prove theorem 1.1 for hypersurfaces of general type, i.e.
the case $n+1-h<0$.  So we let \begin{equation}
n+1+\delta=h
\end{equation}
where integer $\delta\geq 0$.  
\bigskip

Let
\begin{equation}\begin{array}{ccc}
Pr: \mathbf P^{n+\delta} &\dashrightarrow & \mathbf P^n
\end{array}\end{equation}
be the projection. Then
$Pr$ induces a rational map
\begin{equation}\begin{array}{ccc}
\eta: \Gamma_{F_0} &\dashrightarrow & \Gamma_{f_0}\\
c_0^\delta & \rightarrow & Pr\circ c_0.
\end{array}\end{equation}

Let
$f_0\in H^0(\mathcal O_{\mathbf P^n}(h))$ be generic for the $\mathbf P^n$ in (2.133).
Let \begin{equation} 
F_0=f_0+g_0\in H^0(\mathcal O_{\mathbf P^{n+\sigma}}(h))
\end{equation} be generic also.  Let
\begin{equation}
X_0=div(f_0), Y_0=div(F_0).
\end{equation}

We may assume that $c_0^\delta$ lies in the regular locus of $\eta$, and $\eta (c_0^\delta)$ is not a constant map for generic 
$c_0^\delta\in \Gamma_{F_0} $. The map is clearly surjective. Then
if $c_0\in   \Gamma_{f_0}$ is generic, there is an inverse $c_0^\delta=\eta^{-1}(c_0)$ 
in $\Gamma_{F_0}$ which is also generic in $Hom(\mathbf P^1, Y_0)$.
By the result of section 2.1, \begin{equation}
H^1(N_{c_0^\delta/Y_0})=0.
\end{equation}

Let $K$ be the kernel bundle of
the differential 
\begin{equation}\begin{array}{ccc}
Pr_\ast: T_{\mathbf P^{n+\delta}}|_R & \rightarrow T_{\mathbf P^n},
\end{array}\end{equation}
where $R$ is the open set of $\mathbf P^{n+\delta}$ such that
$Pr$ restricted to $R$ is regular. So $K$ is a subbundle of $T_R$.

Now we have an exact sequence of bundles

 \begin{equation}\begin{array}{ccccccccc}
0 &\rightarrow & (c_0^\delta)^\ast (K)  &\rightarrow &  N_{c_0^\delta/ Y_0}&\rightarrow &
N_{c_0/X_0} &\rightarrow & 0.\end{array}\end{equation}
Notice all bundles are over $\mathbf P^1$. 
Therefore we have the exact sequence of vector spaces
\begin{equation}\begin{array}{ccccc}
  H^1( N_{c_0^\delta/ Y_0})&\rightarrow &
H^1 (N_{c_0/X_0}) &\rightarrow & H^2 ( (c_0^\delta)^\ast (K))= 0.\end{array}\end{equation}

By (2.137), $$ H^1( N_{c_0^\delta/ Y_0})=0.$$
Hence $$H^1 (N_{c_0/X_0}) =0.$$

This completes the proof of theorem 1.1.

\bigskip

\vfill\eject

\section{Application 1: First encounter of Voisin's \\
--Calabi-Yau hypersurfaces}

\subsection{Intuitive dimension count} \quad\par
 Our theorem rigorously proves that the following well-known dimension count is valid  except for lower dimensions including K-3 surfaces. 
Consider all rational curves $C$ in a projective space of dimension $n$. 
There are $(n+1)(d+1)-1$ parameters in parameterizing rational curves of degree $d$ in $\mathbf P^n$. For a curve 
to be on a hypersurface of degree $h$, there will be $hd+1$ polynomial equation constraints on these parameters. 
Considering the $3$-dimensional automorphisms group of $\mathbf P^1$, to obtain a rational curve of degree $d$ on the
hypersurface, naively we need to have
\begin{equation} (n+1)(d+1)-1-3\geq hd+1.\end{equation}
For this inequality (3.1),  we expect the strongest on generic hypersurfaces:
 all the mentioned $hd+1$ equations  on the space of rational maps of dimension $(n+1)(d+1)-1$ are algebraically independent, i.e.
they define a reduced variety of codimension $hd+1$.  
Our theorem 1.1 is just a rigorous proof of this expectation for $n\geq 4$. However our proof failed in  the case
$n\leq 3$ because of the vanishing differential form $\omega$. This leads to a fact that the intuitive inequality (3.1) does not hold for $n=3$ (the case of K-3 surfaces). The failure  has no hint in the intuitive observation, but it is confirmed by Chen's result on rational curves on K-3 surfaces [1]. This exceptional case
 for $n=3$ distinguishes itself from other hypersurfaces because of their complex structures,  and therefore it is in its own category.

\bigskip

\bigskip

\subsection{Voisin's conjecture on a Lang's conjecture}\quad\smallskip

In this section, we prove a Vosin's conjecture ([9], p. 114, after conjecture 3.22) which directly disproves a Lang's conjecture on the Abelian covering of
a generic hypersurface. Voisin did not formally list it as a conjecture because this is a generalization of Clemens' conjecture.\par
 (a) Vosin's conjecture (generalization of Clemens' conjecture): in a general Calabi-Yau hypersurface of dimension $\geq 3$, rational 
curves cover a countable union of Zariski closed proper algebraic subsets of codimension $\geq 2$. \par
(b) Lang's  conjecture: if the general hypersurface $X_0$ is not of general type, the union of the images of non-constant rational maps: $A\to  X$
from an Abelian variety $A$ to $X_0$ cover $X_0$.
\bigskip

In corollary 3.33, [9], Voisin proves that if (a) is correct, (b) is incorrect.
 In the following we prove (a) is correct. Therefore Lang's conjecture (b) is incorrect.
\bigskip

\begin{corollary} Let $X_0$ be a Calabi-Yau general hypersurface of $\mathbf P^n, n\geq 4$. 
Let  $m=dim(X_0)-3$.
Then  the dimension of a  parameter space of a family 
of rational curves on $X_0$  can only be $m$ at most.  
\end{corollary} 
\bigskip

\begin{proof} Let $C_0$ be a rational curve and $c_0$ be its normalization. We may assume that
$c_0\in Hom(\mathbf P^1, X_0)$ is generic. 
By theorem 1.1,  $H^1(N_{c_0/X_0})=0$ for generic $c_0\in Hom(\mathbf P^1, X_0)$. 
 Since $$H^1(c_0^\ast(T_{\mathbf P^1}))=0, $$
$$H^1(c_0^\ast(T_{X_0}))=0.$$ 
By theorem 1.2, II.1 in [5],   
 \begin{equation} dim_{[c_0]}Hom(\mathbf P^1, X_0)=-K_{X_0}\cdot C_0+ dim(X_0).\end{equation}
Since $X_0$ is a Calabi-Yau, $c_1(K_{X_0})=0$. Hence $K_{X_0}\cdot C_0=0$. 
We obtain that 
$$dim_{[c_0]}Hom(\mathbf P^1, X_0)= dim(X_0).$$
Since there is a 3-dimensional  group of automorphisms on $\mathbf P^1$, our corollary is proved.
\end{proof}
\par
This immediately implies  the conjecture (a), i.e.  
\bigskip

\begin{corollary}
Rational rational curves on a  general  Calabi-Yau hypersurface of dimension $\geq 3$ cover a countable union of
Zariski closed proper algebraic subsets of codimension $\geq 2$.
\end{corollary}

\bigskip

  \begin{proof}
Because $$dim_{[c_0]}Hom(\mathbf P^1, X_0)=dim(X_0),$$
the dimension of the parameter space of a family of rational curves has to be $$dim(X_0)-3.$$ 
Then each irreducible uniruled variety must have dimension at most $$dim(X_0)-2.$$ 
The corollary is proved. 
\end{proof}

\bigskip

{\bf Remark}.  In the case of $dim(X_0)=3$, corollary 3.2 is the Clemens' conjecture which is proved in [10].
\bigskip

\begin{corollary}
For each $d\geq 1$, the generic hypersurface  $X_0$ of degree $2n-3$ in $\mathbf P^n, n\geq 4$, contains at most 
a finitely number of rational curves of degree $d$. 
\end{corollary}
\bigskip

\begin{proof}
Using the formula (3.2), we obtain that

\begin{equation}\begin{array}{c}  dim_{[c_0]}Hom(\mathbf P^1, X_0)-3= -K_{X_0}\cdot C_0+ dim(X_0)-3. \\
=  (n-4)(1-d) \leq 0.
\end{array}\end{equation}
Considering the 3-dimensional  group of automorphisms on $\mathbf P^1$, our corollary is proved.

\end{proof}

\bigskip

{\bf Remark} The corollary for the case $n\geq 5$ is also a consequence of Voisin's theorem in [7] (see [6]), even though her approach is
quite different from ours.  The corollary for $n=4$ is one of Calabi-Yau cases above, which also is the Clemens' conjecture.

\vfill\eject

\section{Application 2: Second encounter of Vosin's \\
--hypersurfaces of general type}\quad\smallskip

\subsection{ Work by Clemens, Voisin and Pacienza}  \quad\par

Clemens, Voisin, Pacienza and others studied questions about genus of subvarieties of hypersurfaces.  Consequently they obtain many results of rational curves 
on hypersurfaces. 
Theorem 1.1 clearly is different from their general questions on genus of subvarieties, but it has many corollaries that not only coincide with their inequalities but also  go beyond  them.

\bigskip

We first give a bound of degree of hypersurfaces.\bigskip

\begin{corollary} 
Let $X_0$ be a general hypersurface of $\mathbf P^n, n\geq 4$  that contains
an irreducible rational curve $C_0$ of degree $d$, then
\begin{equation} deg(X_0)\leq n+1+{n-4\over d}.\end{equation} 
\end{corollary} 

\bigskip

\begin{proof} Let $h=deg(X_0)$.  Let $c_0: \mathbf P^1\to C_0$ be a birational map (normalization of $C_0$). 
We may assume that
$c_0\in Hom(\mathbf P^1, X_0)$ is generic. 
Because $n\geq 4$, 
$$H^1(N_{c_0/X_0})=0.$$
Next we see that a numerical meaning of $H^1(N_{c_0/X_0})=0$ is the inequality (4.1). \par
This follows from the observation:\par
 $c_0^\ast(T_{X_0})$ is a locally free sheaf over $\mathbf P^1$. Thus it has splitting
$$c_0^\ast(T_{X_0})\simeq \mathcal O_{\mathbf P^1}(a_1)\oplus\cdots\oplus O_{\mathbf P^1}(a_{n-1}).$$
Consider the exact sequence
$$\begin{array}{ccccccccc}
0&\rightarrow & T_{\mathbf P^1} &\rightarrow & c_0^\ast(T_{X_0}) &\rightarrow &
N_{c_0/X_0} &\rightarrow &0. \end{array}$$
We then have 
$$\begin{array}{ccccc}
0=H^1(T_{\mathbf P^1}) &\rightarrow & H^1(c_0^\ast(T_{X_0})) &\rightarrow &
H^1(N_{c_0/X_0}).\end{array} $$
Because $H^1(N_{c_0/X_0})=0$, $H^1(c_0^\ast(T_{X_0}))=0$. 
Hence  $$a_i\geq -1, \ for \ all \ i=1, \cdots, n-1.$$
Notice the tangential automorphisms which lie in $ H^0(c_0^\ast(T_{X_0}))$ have two zeros. This means that at least one of $a_i, i=1, \cdots, {n-1}$ 
is larger than or equal to $2$.   
This shows
$$c_1(c_0^\ast(T_{X_0}))\geq -n+4.$$
By the adjunction formula
$$(n+1-h)d=c_1(c_0^\ast(T_{X_0})).$$
Hence 
$$(n+1-h)d\geq -n+4.$$
This is the inequality (4.1).

\end{proof}

\bigskip

Clemens started the early study of rational curves on hypersurfaces in [2]. Later Voisin and Pancienza borrowed Ein's original idea of using adjunction formula and pushed it further to obtain the minimum genus of 
subvaries in hypersurfaces [6], [7]. Their immediate application of it is for rational curves.  They obtained the following results \par

(1) (Clemens and Voisin)   The maximum degree of a generic hypersurface of $\mathbf P^n$ for $n\geq 4$ containing a rational curve is
\begin{equation}
2n-3
\end{equation}
and the maximum degree can be achieved ( see [2], [7]). \par

(2) (Pancienza) A generic hypersurface of $\mathbf P^n$ of degree $2n-3$ for $n\geq 6$ does not contain rational curves of degree
 larger than 1 (see [6]). 

\bigskip

Our inequality confirms their results and indicates that their numerical results are originated from
\begin{equation}
H^1(N_{c_0/X_0})=0.
\end{equation} 

This is  because if $d\geq 1$, then for all $n\geq 4$ the inequality (4.1) reduces to
\begin{equation}  deg(X_0)\leq n+1+n-4=2n-3.\end{equation}
 This is the maximum degree given by Voisin.  In case  $h=2n-3$ and if $n$ is strictly larger than 3, 
the inequality (4.1) becomes 
$d\leq 1$. This means that such a general hypersurface ($n\geq 5$) does not contain rational curves other than lines.
In the case $n\geq 6$, this is a result (2) by Pacienza. Our result shows that this is also true
for $n=5$. 
\bigskip

The inequality (4.1) is the intuitive inequality (3.1).

\subsection{ Sharp bounds of degrees of rational curves on general hypersurfaces of general type}\quad\par

 It is clear that when $X_0$ is Fano or Calabi-Yau, i.e., $h\leq n+1$, the inequality (4.1) does not produce anything new. 
But the situation is well understood.  If $X_0$ is Fano, it is rationally connected. If it is Calabi-Yau, the situation is discussed in section 3.2.. 
The only remaining case is when $X_0$ is of general type. In this case, there is
a conjecture by Voisin ( conjecture 3.9, [9]):  if a generic hypersurface  is of general type,  
the degrees of rational curves on it are bounded. Our inequality (4.1) shows that 
the conjecture is true, and furthermore we give the sharp bound:
\bigskip

\begin{corollary} 
If $X_0$ is a general hypersurface of $\mathbf P^n, n\geq 4$ of general type, i.e.,
$h> n+1$,  then the degrees of rational curves $C_0$ have an upper bound
\begin{equation} deg(C_0)\leq {n-4\over h-n-1}.\end{equation} 
\end{corollary}

\bigskip
\begin{proof} If $h > n+1$ ( i.e. it is of general type),  the inequality  (4.5)  is exactly the inequality (4.1).  \end{proof}

\bigskip

Above inequalities (4.1) and (4.5) describe  a complete picture of  rational curves on a  general hypersurface $X_0$ of general type in
$\mathbf P^n, n\geq 4$.
We would like to describe them separately because historically some of them were not obtained as one unified result 
from $H^1(N_{c_0/X_0})=0$, but rather as individual results
with other methods. \par
 
(1) if $h$ is in the interval $(2n-3, +\infty)$, the inequality (4.1) for $d\geq 1$ implies that 
 there are no rational curves on $X_0$. This is
the known result of Clemens ([2]) and Voisin ([7], [8]) ;\par
(2) if $h=2n-3$, and $n\geq 5$ the inequality (4.1)   implies that $X_0$ contains no rational curves other than lines.  Pacienza in ([6]) obtained
this result for $n\geq 6$ ;\par
(3) if $h$ is in $({3n\over 2}-1, 2n-3)$, the inequality (4.5) implies that $X_0$ does not contain rational curves other than lines. 
This is the known result of Clemens and Ran ([3]), which is implied by their bound of twisted genus;\par
(4) if $h={3n\over 2}-1$, the inequality (4.5) implies that $X_0$ contains no rational curves other than lines and quadratic curves. This is our new result.\par
(5) if $h$ is in $(n+1, {3n\over 2}-1)$, the degrees of the rational curves $C_0$ are bounded above as in the 
formula (4.5).  This is our new result.

\end{document}